
\parskip=1ex minus .3ex \parindent=0pt
\hoffset=0pt \voffset=0pt
\hsize=149mm \vsize=245mm
\def\ifundef#1{\expandafter\ifx\csname#1\endcsname\relax}
\newif\ifpdftex
\ifundef{pdfoutput}\pdftexfalse\else\pdftextrue\fi
\ifpdftex
 \pdfoutput=1 \pdfcompresslevel=9 \pdfadjustspacing=1
 \pdfpagewidth=210mm \pdfpageheight=297mm 
\fi
\catcode`@=11
\def\usefont[#1]{\font\tmp@fnt = #1\relax\tmp@fnt}
\font\eightrm = cmr10 scaled 800    
\font\sixrm = cmr7 scaled 850
\font\fiverm = cmr5
\font\eightbf = cmbx10 scaled 800
\font\sixbf = cmbx7 scaled 850
\font\fivebf = cmbx5
\font\eightit = cmti10 scaled 800
\font\eighti = cmmi10 scaled 800
\font\sixi = cmmi7 scaled 850
\font\fivei = cmmi5
\font\eightsy = cmsy10 scaled 800
\font\sixsy = cmsy7 scaled 850
\font\fivesy = cmsy5
\font\eightsl = cmsl10 scaled 800
\font\eighttt = cmtt10 scaled 800
\font\fivefk = eufm5                
\font\sixfk = eufm7 scaled 850
\font\sevenfk = eufm7
\font\eightfk = eufm10 scaled 800
\font\tenfk = eufm10
\newfam\fkfam
   \textfont\fkfam=\tenfk \scriptfont\fkfam=\sevenfk
      \scriptscriptfont\fkfam=\fivefk
\def\fk{\fam\fkfam}
\font\fivebbm  = bbmsl10 scaled 500 
\font\sixbbm   = bbmsl10 scaled 600
\font\sevenbbm = bbmsl10 scaled 700
\font\eightbbm = bbmsl10 scaled 800
\font\tenbbm   = bbmsl10
\newfam\bbmfam
   \textfont\bbmfam=\tenbbm \scriptfont\bbmfam=\sevenbbm
      \scriptscriptfont\bbmfam=\fivebbm
\def\bm{\fam\bbmfam}
\def\eightpoint{%
   \textfont0=\eightrm \scriptfont0=\sixrm \scriptscriptfont0=\fiverm
      \def\rm{\fam0\eightrm}%
   \textfont1=\eighti  \scriptfont1=\sixi  \scriptscriptfont1=\fivei
      \def\oldstyle{\fam1\eighti}%
   \textfont2=\eightsy \scriptfont2=\sixsy \scriptscriptfont2=\fivesy
   \textfont\itfam=\eightit \def\it{\fam\itfam\eightit}%
   \textfont\slfam=\eightsl \def\sl{\fam\slfam\eightsl}%
   \textfont\ttfam=\eighttt \def\tt{\fam\ttfam\eighttt}%
   \textfont\bffam=\eightbf \scriptfont\bffam=\sixbf
      \scriptscriptfont\bffam=\fivebf \def\bf{\fam\bffam\eightbf}%
   \textfont\fkfam=\eightfk \scriptfont\fkfam=\sixfk
      \scriptscriptfont\fkfam=\fivefk
   \textfont\bbmfam=\eightbbm \scriptfont\bbmfam=\sixbbm
      \scriptscriptfont\bbmfam=\fivebbm
   \rm}
\skewchar\eighti='177\skewchar\sixi='177
\skewchar\eightsy='60\skewchar\sixsy='60
\def\small{\eightpoint\baselineskip=9.6pt}%
\outer\def\title #1\by #2.{%
 \ifpdftex\pdfinfo{/Title(#1) /Author(#2) /Date(\today)}\fi%
 \centerline{{\usefont[cmbx10 scaled 1440]#1}}%
 \centerline{{\usefont[cmr7](#2)}}%
 \vskip 1ex plus .5ex}
\newcount\secNo \secNo=0
\outer\def\section #1.{%
 \global\advance\secNo by 1
 \global\thmNo=0
 \global\eqnNo=0
 \goodbreak\vskip 3ex\noindent
 {\usefont[cmbx10 scaled 1200]\the\secNo.~#1}
 \vglue 1ex}
\newcount\thmNo \thmNo=0
\def\thmno{\global\advance\thmNo by 1\relax
 \the\secNo.\the\thmNo\space}
\def\today{\number\day\space\ifcase\month\or
 January\or February\or March\or April\or
 May\or June\or July\or August\or
 September\or October\or November\or December\fi
 \space\number\year}
\def\\{\hfill\break}
\let\em\it
\def\htag[#1]#2{\ifpdftex\pdfdest name {#1} xyz\fi\relax#2}
\def\href[#1]#2{\leavevmode\ifpdftex\pdfstartlink
 attr {/Border [0 0 0]} goto name {#1}\relax#2\pdfendlink\else#2\fi}
\everydisplay={\textstyle}
\let\plaineqno\eqno
\newcount\eqnNo \eqnNo=0
\def\eqn@no{(\the\secNo.\the\eqnNo)}
\def\eqno#1$${\global\advance\eqnNo by 1
 \plaineqno\htag[eqn.#1]{{\rm\eqn@no}}$$%
 \ifundef{EQN@#1}%
  \expandafter\xdef\csname EQN@#1\endcsname{\eqn@no}%
 \else
  \errmessage{duplicated equation label}%
 \fi\ignorespaces}
\def\reqn#1{\href[eqn.#1]{{\rm\csname EQN@#1\endcsname}}}
\def\item@exec(#1)(#2,#3,#4)#5{%
 {\parindent=#2\relax\par\leavevmode}%
 \def\argA{#1}\relax\def\argB{#2}\relax%
 \ifx\argA\argB\hangindent=#2\else\hangindent=#3\fi%
 \hangafter=1\llap{#5\enspace}\ignorespaces}%
\def\item@optA[#1]#2{\item@exec(#1)(#1,,){#2}}%
\def\item@optB#1{\item@optA[20pt]{#1}}%
\def\item@opt{%
 \ifx\item@char[\let\item@cmd\item@optA\relax\else%
 \let\item@cmd\item@optB\fi\relax\item@cmd}%
\def\item{\futurelet\item@char\item@opt}

\newcount\foot@no \foot@no=0
\def\footnote{\global\advance\foot@no by 1
  \edef\@sf{\spacefactor\the\spacefactor}%
  \unskip{\raise.4\baselineskip\hbox{
     \usefont[cmr7]\number\foot@no)}}\@sf
  \insert\footins\bgroup
     \everydisplay={}\everypar={}\let\par\endgraf
     \parindent=16pt \parskip=0pt \leftskip=0pt \rightskip=0pt
     \splittopskip=10pt plus 1pt minus 1pt
     \floatingpenalty=20000
     \ifundef{small}\else
        \abovedisplayskip=6pt \abovedisplayshortskip=-4pt
        \belowdisplayskip=4pt \belowdisplayshortskip=4pt
        \small\fi
     \smallskip
     \textindent{\number\foot@no)}
     \bgroup\aftergroup\@foot\let\next}
\def\@foot{\egroup}
\skip\footins=12pt plus 2pt minus 4pt
\dimen\footins=.5\vsize
\catcode`@=12
\def\Z{{\bm Z}}
\def\R{{\bm R}}
\def\C{{\bm C}}
\def\proj(#1){{\bm P}(#1)}
\def\Gl{{\rm Gl}}
\def\SO{{\rm SO}}
\def\gK{K_{int}}

\title Discrete linear Weingarten surfaces
\by F Burstall, U Hertrich-Jeromin \& W Rossman.
\vskip 2em\hbox to\hsize{
 \hfil\vbox{\hsize=0.8\hsize{\small{\bf Abstract.}\enspace
Discrete linear Weingarten surfaces in space forms
are characterized as special discrete $\Omega$-nets,
a discrete analogue of Demoulin's $\Omega$-surfaces.
It is shown that the Lie-geometric deformation of
$\Omega$-nets descends to a Lawson transformation
for discrete linear Weingarten surfaces,
which coincides with the well-known Lawson
correspondence in the constant mean curvature case.
\par}}\hfil}\vskip 1em
\hbox to\hsize{
 \hfil\vbox{\hsize=0.8\hsize{\small{\bf MSC 2010.}\enspace
{\it 53A10\/}, {\it 53C42\/}, 53A40, 37K35, 37K25
\par}}\hfil}\vskip 1em\hbox to\hsize{
 \hfil\vbox{\hsize=0.8\hsize{\small{\bf Keywords.}\enspace
linear Weingarten surface; hyperbolic space; de Sitter space;
Lie sphere geometry; Legendre map; Demoulin's surface;
$\Omega$-surface; deformation; isothermic surface;
Calapso transformation; Lawson correspondence;
discrete isothermic net; mixed area;
constant mean curvature; constant Gauss curvature.
\par}}\hfil}\vskip 1em

\section Introduction.
A Lie geometric approach to flat fronts in hyperbolic space
and, more generally, (smooth) linear Weingarten surfaces in
(Riemannian and Lorentzian) space forms was outlined in two
short notes \href[ref.bjr10]{[7]} and \href[ref.bjr11]{[8]}.
Apart from providing a unified treatment and a natural realm
for a transparent analysis of the singularities of fronts, this
Lie geometric approach also revealed a close relationship to
the theory of isothermic surfaces: linear Weingarten surfaces
in space forms are Lie-applicable,%
\footnote{In fact, linear Weingarten surfaces in Riemannian
  space forms are special Guichard surfaces.}
cf \href[ref.bl29]{[1, \S85]} or \href[ref.muni06]{[16]}.
In particular, non-tubular linear Weingarten surfaces
envelop a pair of isothermic sphere congruences that
separate the curvature sphere congruences harmonically,
see [\href[ref.de11a]{10}; \href[ref.de11b]{11}] and
\href[ref.bl29]{[1, \S85]}, where each isothermic sphere
congruence takes values in a linear sphere complex.  Up to
a mild genericity assumption this yields a characterization
of linear Weingarten surfaces or, more generally, fronts,
see \href[ref.bjr11]{[8]}.

As a main result of the present text we shall provide a similar
characterization in the discrete case,
see \href[thm.lwo]{Thm 2.8},
where discrete ``linear Weingarten nets'' are defined in terms
of mixed areas,
see Defs \href[def.macurvatures]{2.3} and \href[def.lw]{2.4}.
This generalizes and unifies the rather different approaches
to constant mean curvature nets
of \href[ref.bjrs08]{[6]}
and \href[ref.bjl11]{[4]}.

In the process,
see \href[def.omega]{Def 3.1},
we introduce discrete $\Omega$-nets as a new class of
integrable discrete surfaces.
For this definition, we employ the new and geometrically
somewhat obscure idea of K\"onigs duality of suitable
{\em homogeneous coordinates\/}
of K\"onigs nets in a projective space:
classically, {\em affine projections\/} of K\"onigs nets
or surfaces in projective geometry admit duals
--- and are characterized by their existence,
cf \href[ref.bosu08]{[2, Def 2.22]}
or \href[ref.bjl11]{[4, Sect 2]}.
However, this idea is motivated by observations in the
smooth case.

A fact that sets our notion of discrete $\Omega$-nets apart
from their smooth analogues is the existence of multiple
pairs of enveloped isothermic sphere congruences, see
\href[thm.ambiguity]{Lemma 3.3}.
This hints strongly at the non-existence of a sensible
notion of vertex curvature spheres for a discrete Legendre
map (\href[def.legendre]{Def~2.1}) or ``principal contact
element net'' \href[ref.bosu08]{[2, Def 3.23]}:
as initial spheres for a pair of isothermic sphere congruences
of an $\Omega$-net can be chosen arbitrarily in one contact
element, no pair of geometrically defined sphere congruences
will satisfy the aforementioned property of harmonic
separation.  Nevertheless, all isothermic sphere congruences
of the family given in \href[thm.ambiguity]{Lemma 3.3} are
conformal in the sense that they share the same cross ratio
function on faces --- in the smooth case conformality of
the induced metrics is intimately related to the harmonic
separation property.

One merit of describing linear Weingarten surfaces or nets
in the Lie geometric realm is the natural description of their
transformations in terms of the transformations of their
Legendre lifts:
$\Omega$-nets come with their Lie geometric deformation,
the Calapso deformation of \href[def.calapso]{Def 3.9},
as well as with Darboux transformations,
inherited by the corresponding transformations of the enveloped
isothermic sphere congruences.
These transformations give rise to Lawson transformations,
see \href[def.lawson]{Def 4.1},
and Bianchi-B\"acklund transformations of
linear Weingarten nets.

The Lawson transformation will be discussed in detail
in Sect~4 of the present text.
In particular, we shall justify our terminology by showing
that the Lawson transformation becomes the well-known
Lawson correspondence in the case of constant mean
curvature nets,
see \href[expl.cmc]{Example 4.2}
and \href[ref.bjrs08]{[6, Sect 5]}.
The Bianchi-B\"acklund transformation shall be discussed
in a forthcoming paper.

{\it Acknowledgements.\/}
We would like to thank
A Bobenko,
T Hoffmann
and
I Lukyanenko
for fruitful and enjoyable discussions.

\section Discrete linear Weingarten surfaces in space forms.

We aim to describe discrete linear Weingarten surfaces,
defined in terms of mixed areas,
cf \href[ref.bpw10]{[3, Def~8]}
and \href[ref.bjl11]{[4, Def 3.1]},
in Riemannian and Lorentzian space forms in a unified manner.
To this end we consider the space form geometries as
subgeometries of Lie sphere geometry:
fix orthogonal vectors ${\fk p,q}\in\R^{4,2}\setminus\{0\}$
and let
$$
{\fk Q}^3 := \{{\fk y}\in\R^{4,2}\,|\,
{\fk(yy)}=0,{\fk(yq)}=-1,{\fk(yp)}=0\}
\eqno spaceform$$
where $(\cdot\cdot)$ denotes the inner product of $\R^{4,2}$;
$\langle\cdot,\cdots,\cdot\rangle$ will denote the linear span
of vectors.
If $({\fk pp})\neq0$ then ${\fk Q}^3$ is a $3$-dimensional
quadric of constant sectional curvature $-({\fk qq})$.

In this setting, the projective light cone or {\em Lie
quadric\/} ${\cal L}^4 : = \{ \langle {\fk y}\rangle\,|\,
{\fk y}\in\R^{4,2},({\fk y}{\fk y})=0 \}\subset\proj(\R^{4,2})$ parametrises
the set of oriented $2$-spheres (thus, complete, totally
umbilic hypersurfaces) in ${\fk Q}^3$ via
$$
s\mapsto {\fk Q}^3\cap s^\perp.
$$
In particular, for ${\fk y}\in{\fk Q}^3$, $s=\langle {\fk y}\rangle$
corresponds to the point-sphere $\{{\fk y}\}$ while, when
$s\in{\cal L}^4$ differs from its reflection $s'$ in the
hyperplane orthogonal to ${\fk p}$, $s,s'$ correspond to the
same sphere but with opposite orientations.

In general, a non-zero point ${\fk k}\in\R^{4,2}$ (or,
more properly, a point $\langle{\fk
k}\rangle\in\proj(\R^{4,2})$), defines the {\em linear
sphere complex\/} ${\cal L}^4\cap {\fk k}^\perp$, a
$3$-dimensional family of $2$-spheres.  In particular,
taking ${\fk k}={\fk q}$ yields
$$
{\fk P}^3 := \{{\fk y}\in\R^{4,2}\,|\,
{\fk(yy)}=0,{\fk(yq)}=0,{\fk(yp)}=-1\},\eqno hyperplane
$$
the complex of (spacelike if $({\fk pp})>0$) hyperplanes
(thus, complete, totally geodesic hypersurfaces) in
the space form ${\fk Q}^3$, cf \href[ref.imdg]{[14, Sect 1.4]}.

Two oriented $2$-spheres are in oriented contact if and only
if the corresponding points of ${\cal L}^4$ are orthogonal.
It follows that lines in ${\cal L}^4$ correspond to pencils
of $2$-spheres sharing a common contact element and so
parametrise those contact elements.  For more details, see
Cecil \href[ref.cecil]{[9, Chapter~1]}. 

To make this approach more tangible,
assume that $({\fk pp})=\mp 1$ and $({\fk qq})\neq 0$.
Now the constant offset
$$
  \{ {\fk x=f+{q\over(qq)}} \,|\, {\fk f\in Q}^3 \}
  \subset\langle{\fk p,q}\rangle^\perp
$$
yields the standard model of a space form as a
(connected component of a) quadric in a $4$-dimensional
linear space with non-degenerate inner product,
and the unit (timelike if $({\fk pp})=+1$) tangent space
of the space form at ${\fk x=f+{q\over(qq)}}$
becomes the constant offset
$$
  \{ {\fk n=t+{p\over(pp)}} \,|\,
     {\fk t\in P}^3\cap\langle{\fk f}\rangle^\perp \}
  \subset\langle{\fk p,q}\rangle^\perp.
$$
In the case $({\fk qq})=0$ of a flat ambient space form
geometry, the situation becomes slightly less obvious:
here a choice of origin ${\fk o\in Q}^3$ yields an
identification via inverse stereographic projection,
$$
  {\fk Q}^3\ni{\fk f=o+x}+{1\over2}{\fk(xx)q}
    \enspace\leftrightarrow\enspace
  {\fk x}\in\langle{\fk o,p,q}\rangle^\perp
  \cong \cases{
    \R^3 & if $({\fk pp})=-1$, \cr
    \R^{2,1} & if $({\fk pp})=+1$, \cr
  }
\eqno stereoproj$$
and ${\fk P}^3\cap\langle{\fk f}\rangle^\perp$ becomes the
unit (timelike in the Lorentzian case) tangent space of
$\R^3$ or $\R^{2,1}$, respectively,
via
$$
  {\fk P}^3\cap\langle{\fk f}\rangle^\perp \ni
    {\fk t=-{p\over(pp)}+n+(xn)q}
    \enspace\leftrightarrow\enspace
  {\fk n}\in\langle{\fk o,p,q}\rangle^\perp.
\eqno stereoprojt$$

Now consider a discrete principal (circular) net%
\footnote{For simplicity we restrict to $\Z^2$ as a domain;
  throughout, $\Z^2$ may be replaced by a (simply connected)
  quad-graph.}
${\fk f}:\Z^2\to{\fk Q}^3$,
that is, ${\fk f}$ has planar faces in $\R^{4,2}$,
cf \href[ref.bosu08]{[2, Thm 3.9]}.
For non-degeneracy we assume that neither edges nor diagonals
of ${\fk f}$ be isotropic:
if $(ijkl)$ denotes an elementary quadrilateral of $\Z^2$ then
the vectors
$$
  d{\fk f}_{ij} := {\fk f}_j-{\fk f}_i
    \enspace{\rm and}\enspace
  \delta{\fk f}_{ik} := {\fk f}_k-{\fk f}_i
$$
will be assumed to be non-null.
In particular, any two or three vertices of a face of ${\fk f}$
span a $2$- or $3$-dimensional subspace of $\R^{4,2}$,
respectively, with a non-degenerate induced inner product.
Further, in order to be able to define the Gau\ss\ and mean
curvatures via mixed areas below, we shall assume that the
faces of ${\fk f}$ have non-parallel diagonals,
so that their areas do not vanish.

Such a principal net admits a $2$-parameter family of
Gau\ss\ maps, that is, unit (timelike in the Lorentzian case)
``normal'' vector fields along ${\fk f}$ so that,
for each edge $(ij)$, there is an edge curvature sphere
$\kappa_{ij}$ that is orthogonal to the ``normal'' vectors
at the endpoints%
\footnote{Note that, in contrast to the Euclidean case,
  the normal lines in ${\fk Q}^3$ defined by the Gau\ss\ map
  do not necessarily intersect:
  this is the case when the ``curvature sphere'' $\kappa_{ij}$
  is not a distance sphere in the ambient space form geometry.}
${\fk f}_i$ and ${\fk f}_j$,
cf \href[ref.bosu08]{[2, Thm 3.36]}.
In our Lie geometric setup, a choice of Gau\ss\ map for
the principal net ${\fk f}$  amounts to a choice of a
``tangent plane'' congruence ${\fk t}:\Z^2\to{\fk P}^3$
with ${\fk t\perp f}$.
This pair of maps gives rise to the Legendre lift
(principal contact element net,
cf \href[ref.bosu08]{[2, Def 3.23]})
of a principal net in a space form with Gau\ss\ map:

\proclaim\htag[def.legendre]{\thmno Def}.
Let ${\fk f}:\Z^2\to{\fk Q}^3$ be a principal net in a quadric
${\fk Q}^3$ of constant sectional curvature with tangent plane
congruence ${\fk t}:\Z^2\to{\fk P}^3$, ${\fk t\perp f}$.
The line congruence%
\footnote{Here, a ``line congruence'' just means a map
into the space of lines in a projective space.}
$$
\Z^2\ni i\mapsto f_i:=\langle{\fk f}_i,{\fk t}_i\rangle
\subset{\cal L}^4 
$$
will be
called the {\em Legendre lift\/} of the pair $({\fk f,t})$
if adjacent lines $f_i$ and $f_j$ intersect;
$\kappa_{ij}:=f_i\cap f_j$ is called the
{\em curvature sphere\/} of $f$ on the edge $(ij)$.\\
The pair $({\fk f,t})$ will be called the
{\em space form projection\/} of the Legendre map $f$.

We shall exclusively deal with pairs
$({\fk f,t}):\Z^2\to{\fk Q}^3\times{\fk P}^3$
occuring as space form projections of Legendre maps.
Note that, generically,%
\footnote{In the definite case $({\fk pp})<0$,
  the only obstruction is that the point sphere map
  may hit the infinity boundary of the space form;
  in the Lorentzian case, additional obstructions
  occur as a contact element may consist entirely
  of ``point spheres''.}
any choice of a point sphere complex ${\fk p}\in\R^{4,2}$,
$({\fk pp})\neq 0$, and a space form vector
${\fk q}\in\R^{4,2}\setminus\{0\}$, $({\fk qp})=0$,
gives rise to a space form projection $({\fk f,t})$
of a given Legendre map $f$.

As the edge curvature sphere $\kappa_{ij}\in{\cal L}^4$ is
obtained as the intersection $\kappa_{ij}=f_i\cap f_j$ of
the lines of the Legendre lift of a principal net ${\fk f}$
with tangent plane congruence ${\fk t}$ at the endpoints
of an edge $(ij)$, we may write
$$
  \kappa_{ij}
  = {\fk t}_i + k_{ji}{\fk f}_i
  = {\fk t}_j + k_{ij}{\fk f}_j
\eqno cs$$
for (a lift of) the curvature sphere with suitable coefficients
$k_{ij}$ and $k_{ji}$.
Now
$$
  k_{ji}
  = {(\kappa_{ij}{\fk q})\over(\kappa_{ij}{\fk p})}
  = k_{ij},
$$
showing that $(ij)\mapsto k_{ij}$ is an {\em edge function\/},
that is, takes equal values for opposite orientations of
an edge.
This yields a notion of a {\em principal curvature function\/}
on the edges of a principal net ${\fk f}$ in ${\fk Q}^3$
with tangent plane congruence ${\fk t}$.
Rewriting \reqn{cs} as
$$
  0 = d{\fk t}_{ij} + k_{ij}d{\fk f}_{ij}
\eqno rodrigues$$
we obtain a Rodrigues' type formula.
Conversely, \reqn{rodrigues} implies that ${\fk f}$,
hence ${\fk t}$, are conjugate nets in $\R^{4,2}$
as long as the principal curvature function $k$
is not constant around an elementary quadrilateral,
that is, away from {\em umbilical faces\/},
where the curvature spheres of the four edges of
a face coincide.
Thus we obtain the following characterization:

\proclaim\htag[thm.rodrigues]{\thmno Lemma}.
A space form projection
$({\fk f,t}):\Z^2\to{\fk Q}^3\times{\fk P}^3$
of a Legendre map $f$ is a pair of edge-parallel
nets in $\R^{4,2}$.
Conversely, if ${\fk f}:\Z^2\to{\fk Q}^3$ and
${\fk t}:\Z^2\to{\fk P}^3$ satisfy \reqn{rodrigues},
then any non-umbilical face of ${\fk f}$ is planar
and, away from umbilical faces, $({\fk f,t})$ is
the space form projection of a Legendre map $f$.

In particular, the faces of ${\fk f}$ and ${\fk t}$ lie in
parallel planes so that the $\Lambda^2\R^{4,2}$-valued
(mixed) area functions
$$
  A({\fk t,t})_{ijkl}
  = {1\over 2}\delta{\fk t}_{ik}\wedge\delta{\fk t}_{jl}
    \enspace{\rm and}\enspace
  A({\fk f,t})_{ijkl}
  = {1\over 4}\{
      \delta{\fk f}_{ik}\wedge\delta{\fk t}_{jl}
    + \delta{\fk t}_{ik}\wedge\delta{\fk f}_{jl}
    \}
$$
are multiples of
$
  A({\fk f,f})_{ijkl}
  = {1\over 2}\delta{\fk f}_{ik}\wedge\delta{\fk f}_{jl};
$
note that $A({\fk f,f})\neq0$ by our regularity assumption
on ${\fk f}$.

\proclaim\htag[def.macurvatures]{\thmno Lemma \& Def}.
There are two functions, $H$ and $K$, defined on the faces%
\footnote{Note that these Gau\ss\ and mean curvatures do not
  depend on the orientation of an elementary quadrilateral
  $(ijkl)$:
  reversing the orientation, all (mixed) areas change sign
  so that the curvatures remain unaffected.}
of a space form projection of a Legendre map so that
$$
  0 \equiv A({\fk f,t}) + H\,A({\fk f,f})
         = A({\fk t,t}) - K\,A({\fk f,f}).
$$
These will be called the {\em mean curvature\/} and
{\em Gau\ss\ curvature\/} of the pair $({\fk f,t})$,
respectively.

As the mixed areas are invariant under translation,
the mean and Gau\ss\ curvatures defined here clearly
coincide with those of \href[ref.bjl11]{[4, Def 3.1]}
in the case of a Riemannian ambient geometry.

To see that they coincide with the ones of
\href[ref.bosu08]{[2, Def 4.45]} and
\href[ref.bpw10]{[3, Def 8]}
in the case of a principal net ${\fk x}$ in $\R^3$
with (unit) Gau\ss\ map ${\fk n}$ we employ \reqn{stereoproj}
and \reqn{stereoprojt} to observe that the mixed areas of
$$
  {\fk f} = {\fk o+x}+{1\over 2}{\fk(xx)q}
    \enspace{\rm and}\enspace
  {\fk t} = {\fk p+n+(xn)q}
$$
take values in
$
  \Lambda^2\langle{\fk o,p,q}\rangle^\perp
  \oplus
  (\langle{\fk o,p,q}\rangle^\perp\wedge\langle{\fk q}\rangle)
$
and that the mean and Gau\ss\ curvatures $H$ and $K$ are
therefore determined by the
$\Lambda^2\langle{\fk o,p,q}\rangle^\perp$-parts
$A({\fk n,n})$, $A({\fk x,n})$ and $A({\fk x,x})$
of the mixed areas
$A({\fk t,t})$, $A({\fk t,f})$ and $A({\fk f,f})$.

\proclaim\htag[def.lw]{\thmno Def}.
The space form projection
$({\fk f,t}):\Z^2\to{\fk Q}^3\times{\fk P}^3$
of a Legendre map is called a {\em linear Weingarten net\/}
if its mean and Gau\ss\ curvatures satisfy a non-trivial
affine relation
$$
  0 = \alpha\,K + 2\beta\,H + \gamma.
\eqno lw$$

Note the symmetry of the situation:
in case $({\fk qq})\neq 0$ we may interchange the geometric
interpretations of ${\fk p}$ and ${\fk q}$,
thus swapping the roles of ${\fk f}$ and ${\fk t}$.
That is, ${\fk t}$ is interpreted as the principal net
and ${\fk f}$ as its tangent plane congruence.
As long as $A({\fk t,t})\neq 0$, that is, $K\neq 0$,
we obtain ${H\over K}$ and ${1\over K}$ as the mean
and Gau\ss\ curvatures of the pair $({\fk t,f})$,
which is therefore a linear Weingarten net also.
In particular, if $({\fk f,t})$ is a minimal net, $H\equiv 0$,
then so is $({\fk t,f})$.
In this case $A({\fk f,t})\equiv 0$ so that
${\fk f}$ and ${\fk t}$ are K\"onigs dual nets in $\R^{4,2}$,
cf \href[ref.bosu08]{[2, Def 2.22]}:

\proclaim\htag[def.dual]{\thmno Def}.
Two discrete maps $\sigma^\pm$ into an affine space are
called {\em K\"onigs dual\/} if they are edge-parallel
and their opposite diagonals are parallel:
for any edge $(ij)$ and any elementary quadrilateral $(ijkl)$
$$
  d\sigma^+_{ij}\parallel d\sigma^-_{ij}
    \enspace{\sl and}\enspace
  \delta\sigma^\pm_{ik}\parallel\delta\sigma^\mp_{jl}.
$$

Indeed,
$
    \delta\sigma^+_{ik}\wedge\delta\sigma^-_{jl}
  = \delta\sigma^-_{ik}\wedge\delta\sigma^+_{jl}
$
for any edge-parallel nets $\sigma^+$ and $\sigma^-$,
so that the vanishing of their mixed area,
$
  A(\sigma^+,\sigma^-)_{ijkl} \equiv 0,
$
is readily seen to be equivalent to their opposite
diagonals being parallel.
Thus we obtain a characterization of minimal nets in space
forms via K\"onigs duality,
cf \href[ref.bjl11]{[4]}:

\proclaim\htag[thm.minimal]{\thmno Thm}.
The space form projection
$({\fk f,t}):\Z^2\to{\fk Q}^3\times{\fk P}^3$
of a Legendre map is minimal if and only if
${\fk f}$ and ${\fk t}$ are K\"onigs dual lifts
in $\R^{4,2}$ of nets in the Lie quadric.

We aim to generalize this description for
linear Weingarten nets.
To this end, suppose that $\sigma^\pm$ are K\"onigs dual lifts
of sphere congruences\footnote{A sphere
congruence is simply a map $s:\Z^2\to{\cal L}^4$ to the
space of $2$-spheres: the components ${\fk f,t}$ of a
space-form projection of a Legendre map are examples.}
$s^\pm$ spanning a Legendre map $f$;
further suppose that each sphere congruence takes values in
a linear sphere complex ${\fk k}^\pm$.
Since $\sigma^+$ and $\sigma^-$ are edge-parallel nets
and $\sigma^\pm\perp{\fk k}^\pm$, the inner products
$(\sigma^\pm{\fk k}^\mp)\equiv const$.
As long as these inner products do not vanish we may,
without loss of generality, assume the same relative
normalizations as for space form projections:%
\footnote{We shall see below that we need to allow $\sigma^\pm$
  to become complex conjugate in order to capture general
  linear Weingarten nets:
  in this case ${\fk k}^\pm$ can be assumed to be complex
  conjugate as well and this relative normalization can
  be achieved while maintaining complex conjugacy since
  $(\sigma^\pm{\fk k}^\mp)$ are complex conjugate also.}
$$
  (\sigma^\pm{\fk k}^\pm) = 0
    \enspace{\rm and}\enspace
  (\sigma^\pm{\fk k}^\mp) = -1.
$$

Now ${\fk k}^\pm$ span a plane and choosing
a point sphere complex ${\fk p}$ and
a space form vector ${\fk q}$
for a space form projection $({\fk f,t})$ of $f$ in this plane,
$\langle{\fk q,p}\rangle=\langle{\fk k^+,k^-}\rangle$,
our relative normalizations control the relation between
basis transformations:
with $B\in\Gl(2)$ a change of basis
$$
  ({\fk q,p}) = ({\fk k^+,k^-})B
    \enspace{\rm yields}\enspace
  (\sigma^-,\sigma^+) = ({\fk f,t})B^t.
\eqno cbas$$

As both the symmetric products on $\langle{\fk k^-,k^+}\rangle$
and the mixed areas of pairs of edge-parallel nets spanning $f$,
are symmetric bilinear forms they change in a similar way:
$$
  \left(
  {
    {\fk q\odot q}\atop{\fk p\odot q}
  }\,{
    {\fk q\odot p}\atop{\fk p\odot p}
  }
  \right)
  = B^t
  \left(
  {
    {\fk k^+\odot k^+}\atop{\fk k^-\odot k^+}
  }\,{
    {\fk k^+\odot k^-}\atop{\fk k^-\odot k^-}
  }
  \right)
  B
    \enspace{\rm and}\enspace
  \left({
    A(\sigma^-,\sigma^-)\atop A(\sigma^+,\sigma^-)
    }\,{
    A(\sigma^-,\sigma^+)\atop A(\sigma^+,\sigma^+)
  }\right)
  = B
  \left({
    A({\fk f,f})\atop A({\fk t,f})
    }\,{
    A({\fk f,t})\atop A({\fk t,t})
  }\right)
  B^t.
\eqno cdyad$$

Thus if $\sigma^\pm$ are K\"onigs dual,
$A(\sigma^+,\sigma^-)\equiv 0$,
then the constructed space form projection $({\fk f,t})$
is a linear Weingarten net:
$$
  \alpha A({\fk t,t})
  - 2\beta A({\fk t,f})
  + \gamma A({\fk f,f})
  \equiv 0
\eqno malw$$
for suitable constants $\alpha,\beta,\gamma\in\R$
which are determined from the basis representation of
the symmetric bilinear form%
\footnote{In the smooth case, $W$ realizes the linear
  Weingarten condition as an orthogonality condition
  for the curvature spheres.
  Note that, when ${\fk k}^\pm$ are complex conjugate,
  $W$ is real.}
$$
  W := 2\,{\fk k^-\odot k^+}
  = {1\over2(\alpha\gamma-\beta^2)}\{
    \alpha\,{\fk q\odot q}
  + 2\beta\,{\fk q\odot p}
  + \gamma\,{\fk p\odot p}
  \}.
\eqno lwtensor$$

To see the converse we merely reverse this line of argument.
Let $({\fk f,t}):\Z^2\to{\fk Q}^3\times{\fk P}^3$ be a
linear Weingarten net, that is, its mixed areas satisfy
a linear relation \reqn{malw}.
We seek ${\fk k}^\pm$ satisfying \reqn{lwtensor},
that is, factorizing
$$
  W :
  = \alpha\,{\fk q\odot q}
  + 2\beta\,{\fk q\odot p}
  + \gamma\,{\fk p\odot p}.
$$

Clearly, this ambition is in vain if $W$ does not have
full rank, that is, if $\alpha\gamma-\beta^2=0$:
as the sought-after ${\fk k}^\pm$ are linearly independent
${\fk k^+\odot k^-}$ has rank $2$.
Thus we shall exclude this case from the investigation.
The following terminology is chosen in analogy to the
smooth case, where the linear Weingarten surfaces with
$\alpha\gamma=\beta^2$ are those with a constant
principal curvature:

\proclaim\htag[def.tubular]{\thmno Def}.
A linear Weingarten net
$({\fk f,t}):\Z^2\to{\fk Q}^3\times{\fk P}^3$
with $\alpha\gamma-\beta^2=0$ will be called
{\em tubular\/}.

In the non-tubular case $\delta^2:=\beta^2-\alpha\gamma\neq 0$
we can now solve the factorization problem up to order and
(geometrically irrelevant) scaling,
hence obtaining a pair of K\"onigs dual lifts $\sigma^\pm$
of sphere congruences $s^\pm$:

\hangindent=2em\hangafter=1
If $\alpha\neq0$
  then ${\fk k^+\odot k^-}={1\over 4\alpha}W$ with
  ${\fk k}^\pm:
    = {1\over2\alpha}\{{\fk(\alpha q+\beta p)\pm\delta p}\}$
  and $\sigma^\pm={\fk f\pm{1\over\delta}(\beta f-\alpha t)}$
  yield the sought-after K\"onigs dual lifts of sphere
  congruences
  $s^\pm:=\langle\sigma^\pm\rangle:\Z^2\to{\cal L}^4$.
  By construction, the sphere congruences $s^\pm$ take values
  in (different) linear sphere complexes ${\fk k}^\pm$,
  $s^\pm\perp{\fk k}^\pm$.
  Note that ${\fk k}^\pm$ and $\sigma^\pm$ become
  complex conjugate when $\alpha\gamma-\beta^2>0$.

\hangindent=2em\hangafter=1
If $\alpha=0$
  then $\beta\neq0$ and ${\fk k^+\odot k^-}={1\over2\beta}W$
  with ${\fk k}^-:={\fk p}$ and
  ${\fk k}^+:={\fk q}+{\gamma\over2\beta}{\fk p}$.
  In this case, not too surprisingly, we recover the constant
  mean curvature net $\sigma^-={\fk f}$ together with its
  {\em mean curvature sphere congruence\/}
  $\sigma^+={\fk t}+H\,{\fk f}$,
  cf \href[ref.bjrs08]{[6, Def 5.1]}
  or \href[ref.bjl11]{[4, Def 4.1]},
  as a pair of enveloped sphere congruences with
  K\"onigs dual lifts.
  Again, $s^\pm$ take values in the linear sphere complexes
  given by ${\fk k}^\pm$.

Thus we have proved:

\proclaim\htag[thm.lwo]{\thmno Thm}.
The Legendre lift of a non-tubular linear Weingarten net
is spanned by a pair of (possibly complex conjugate)
sphere congruences $s^\pm$ that admit K\"onigs dual lifts.
The sphere congruences take values in different linear
sphere complexes ${\fk k}^\pm$, $s^\pm\perp{\fk k}^\pm$.\\
Conversely, if $f$ is a Legendre map spanned by a pair of
sphere congruences $s^\pm$ that admit K\"onigs dual lifts
$\sigma^\pm$ and take values in different linear sphere
complexes ${\fk k}^\pm$,
then any space form projection $({\fk f,t})$ of $f$ with
$\langle{\fk q,p}\rangle=\langle{\fk k^+,k^-}\rangle$
is a non-tubular linear Weingarten net.

In particular, as parallel nets in a space form are
obtained from space form projections
$({\fk f,t})$ and $(\tilde{\fk f},\tilde{\fk t})$
of the same Legendre map $f$ with respect to bases
$({\fk q,p})$ and $(\tilde{\fk q},\tilde{\fk p})$
that are related by an orthogonal transformation
of their common plane,
we have not too surprisingly also learned:

\proclaim\htag[thm.lwparallel]{\thmno Cor}.
The parallel nets of a linear Weingarten net in a space form
are linear Weingarten.

\section Discrete \char 10-surfaces and their Calapso deformation.

Recall that a discrete {\em Legendre map\/} is a line
congruence
$\Z^2\ni i\mapsto f_i\subset{\cal L}^4\subset\proj(\R^{4,2})$
so that adjacent lines share a (unique)
{\em curvature sphere\/},
$f_i\cap f_j=\kappa_{ij}\in{\cal L}^4$,
cf \href[ref.bosu08]{[2, Def 3.23]}.
In \href[thm.lwo]{Thm 2.8} we have seen that non-tubular
linear Weingarten nets lift to Legendre maps that are
spanned by pairs of sphere congruences admitting
K\"onigs dual lifts,
that is, they lift to $\Omega$-nets of Lie sphere geometry:

\proclaim\htag[def.omega]{\thmno Def}.
A discrete Legendre map is called a
{\em discrete $\Omega$-net\/}
if it is spanned by a pair of sphere congruences
$s^\pm:\Z^2\to{\cal L}^4$ that admit K\"onigs dual lifts
$\sigma^\pm:\Z^2\to\R^{4,2}$.

For regularity we assume, as for the principal net ${\fk f}$
of a space form projection, that the spheres at different
vertices of an elementary quadrilateral do not touch,
that is, $d\sigma^\pm$ and $\delta\sigma^\pm$ never
become isotropic.
Hence
the endpoints of an edge of $s^\pm$
span a $2$-dimensional Minkowski space and
the vertices of any face of $s^\pm$
span a $3$-dimensional space with non-degenerate
induced inner product.
We shall also exclude umbilical faces, where the
curvature spheres of the incident edges all coincide.

As an immediate consequence of this definition,
the sphere congruences $s^\pm$ are
{\em Ribaucour sphere congruences\/}
in the sense of \href[ref.bosu08]{[2, Def 3.27]}:
both sphere congruences $s^\pm$ have planar faces.

Moreover, $s^\pm:\Z^2\to{\cal L}^4$ are
{\em isothermic sphere congruences\/} as K\"onigs nets
in the Lie quadric.%
\footnote{Note our somewhat unusual point of view:
  naturally, K\"onigs nets form a class of nets in
  projective geometry as they can be characterized
  in terms of incidence relations,
  while the notion of their duality belongs to
  an affine subgeometry of the projective ambient
  geometry as it relies on a notion of parallelity.
  In contrast, we consider K\"onigs duality of
  (not necessarily affine) lifts of K\"onigs nets
  in the linear space of homogeneous coordinates
  of their ambient projective geometry.}
To see that $s^\pm$ are indeed K\"onigs nets in
$\proj(\R^{4,2})$ we employ a characterization of
K\"onigs nets in terms of their diagonal vertex stars,%
\footnote{Alternatively, planarity of intersection points
  of diagonals of adjacent faces could be employed,
  cf \href[ref.bosu08]{[2, Thm 2.26]}.}
cf \href[ref.bosu08]{[2, Thm 2.27]}.
Since the K\"onigs dual lifts $\sigma^\pm$ of $s^\pm$ are
K\"onigs nets (in $\R^{4,2})$ their diagonal vertex stars
lie in $3$-dimensional (affine) subspaces of $\R^{4,2}$.
Consequently, the diagonal vertex stars of $s^\pm$ span
$4$-dimensional (linear) subspaces of $\R^{4,2}$,
hence lie in $3$-dimensional (projective) subspaces
of $\proj(\R^{4,2})$. To summarise:

\proclaim\htag[thm.isothermic]{\thmno Lemma}.
If $f=s^-\oplus s^+$ is an $\Omega$-net spanned by
a pair of sphere congruences $s^\pm$ that admit
K\"onigs dual lifts,
then $s^\pm$ are isothermic and, in particular,
Ribaucour sphere congruences.

Below we shall see that the sphere congruences $s^\pm$
come with a cross ratio factorizing edge-labelling
as well as their respective isothermic loops
of flat connections,
cf \href[ref.do06]{[12, Def 4 and Prop 10]}
and \href[ref.bjrs08]{[6, Def 2.1 and Lemma 2.5]}.

The definition of an $\Omega$-net aims to provide a discrete
analogue of smooth $\Omega$-surfaces, the generic%
\footnote{Excluding ``$\Omega_0$-surfaces'',
  where the two isothermic sphere congruences coincide
  with one of the two curvature sphere congruences.}
deformable surfaces of Lie geometry,
see \href[ref.bl29]{[1, \S85]} or \href[ref.muni06]{[16]}.
In the smooth case, these come in two classes:
the ones originally investigated by Demoulin
  [\href[ref.de11a]{10}; \href[ref.de11b]{11}],
  which are given by a real pair of isothermic
  sphere congruences,
and the ones where the enveloped isothermic
  sphere congruences become complex conjugate.
In the discrete case, these two classes merge and
\href[def.omega]{Def 3.1} captures the entire class
without the need to allow for complex conjugate pairs
of sphere congruences:

\proclaim\htag[thm.ambiguity]{\thmno Lemma}.
Let $f$ be an $\Omega$-net.
Then, for any given pair of spanning spheres $s^\pm_0\in f_0$
at an initial point $0\in\Z^2$,
there is a pair of isothermic sphere congruences $s^\pm$
through $s^\pm_0$
that admit K\"onigs dual lifts and span $f$.

To prove this lemma we investigate how to construct a new
pair of isothermic sphere congruences from a given one:
thus let $f=s^+\oplus s^-$ with a pair of isothermic sphere
congruences $s^\pm$ that have K\"onigs dual lifts $\sigma^\pm$.
Hence there is a real function $r$ so that $\sigma^\pm$
satisfy the Christoffel formula
$$
  d\sigma^-_{ij} = r_ir_j\,d\sigma^+_{ij}
    \enspace\Leftrightarrow\enspace
  d\sigma^+_{ij} = {1\over r_ir_j}\,d\sigma^-_{ij}
\eqno christoffel$$
and $\mu^\pm=r^{\pm1}\sigma^\pm$ are Moutard lifts of $s^\pm$,
see \href[ref.bosu08]{[2, Thms 2.31 and 2.32]}:%
\footnote{By the $\pm$-symmetry of \reqn{christoffel}
  the two functions $r^\pm$ obtained
  from \href[ref.bosu08]{[2, Thms 2.31 and 2.32]}
  can be chosen to be reciprocal, $r^\pm=r^{\pm 1}$
  with a single function $r$.}
note that the Moutard equations \href[ref.bosu08]{[2, (2.44)]}
for $\mu^\pm$ are nothing but the integrability conditions%
\footnote{Having excluded umbilical faces we must have
  $(r_k-r_i)(r_l-r_j)\neq0$.}
for \reqn{christoffel},
$$
    (r_k^{\pm1}-r_i^{\pm1})
      \{r_l^{\pm1}\sigma^\pm_l-r_j^{\pm1}\sigma^\pm_j\}
  = (r_l^{\pm1}-r_j^{\pm1})
      \{r_k^{\pm1}\sigma^\pm_k-r_i^{\pm1}\sigma^\pm_i\}.
$$
Now, for any two constants $c^+\neq c^-$,
$$
  \tilde\sigma^\pm :
  = {1\over r+c^\mp}\{\sigma^-+c^\pm r\sigma^+\}
\eqno ambiguity$$
yields K\"onigs dual lifts of another pair of isothermic
sphere congruences spanning the same $\Omega$-net:
the fact that $\tilde\sigma^\pm$ are edge-parallel,
$$
    (r_i+c^-)(r_j+c^-)\,d\tilde\sigma^+_{ij}
  = (r_i+c^+)(r_j+c^+)\,d\tilde\sigma^-_{ij},
$$
hinges on the Christoffel equation \reqn{christoffel},
while the K\"onigs duality of $\tilde\sigma^\pm$,
$$
    {(r_i+c^-)(r_k+c^-)\over r_k-r_i}\,
    \delta\tilde\sigma^+_{ik}
  = {(r_j+c^+)(r_l+c^+)\over r_l-r_j}\,
    \delta\tilde\sigma^-_{jl},
$$
follows from the K\"onigs duality
\href[ref.bosu08]{[2, (2.40)]}
of $\sigma^\pm$ and the Moutard equations
\href[ref.bosu08]{[2, (2.44)]}
for $\mu^\pm$,
$$
  {r^{\pm1}_ir^{\pm1}_k\over r^{\pm1}_k-r^{\pm1}_i}\,
    \delta\sigma^\pm_{ik}
  = {1\over r^{\pm1}_l-r^{\pm1}_j}\,\delta\sigma^\mp_{jl}
      \enspace{\rm and}\enspace
    {1\over r^{\pm1}_k-r^{\pm1}_i}\,\delta\mu^\pm_{ik}
  = {1\over r^{\pm1}_l-r^{\pm1}_j}\,\delta\mu^\pm_{jl}.
\eqno kmeq$$
These computations also yield
$\tilde r:={r+c^-\over r+c^+}$ as a rescaling for
Moutard lifts $\tilde\mu^\pm=\tilde r^{\pm1}\tilde\sigma^\pm$
of the new isothermic sphere congruences
$\tilde s^\pm=\langle\tilde\sigma^\pm\rangle$.

Thus, given an $\Omega$-net $f=s^+\oplus s^-$ in terms of
a pair of isothermic sphere congruences $s^\pm$ with
K\"onigs dual lifts $\sigma^\pm$,
another such pair $\tilde s^\pm$ can be constructed to
pass through any two spanning spheres at a given initial
point $0\in\Z^2$ by choosing the constants $c^\pm$
appropriately.
We have therefore proved Lemma~\href[thm.ambiguity]{3.3}.

Clearly, choosing $c^\pm$ complex conjugate
in \reqn{ambiguity},
say $c^\pm=\pm i$, a complex conjugate pair
$$
  \tilde\sigma^\pm
  = {1\over r\mp i}\{\sigma^-\pm ir\sigma^+\}
$$
of K\"onigs dual lifts in $\R^{4,2}\otimes\C\cong\C^6$
is obtained from a real pair $\sigma^\pm$.
To see that, conversely, a real pair can be obtained from
a complex conjugate pair, $\sigma^+=\overline{\sigma^-}$,
first note that we can, without loss of generality,
assume that $|r|^2\equiv 1$:
on any edge $|r_i|^2|r_j|^2=1$ by \reqn{christoffel}
when $\sigma^\pm$ are complex conjugate;
hence $|r|^2\equiv 1$ as soon as the scaling of $r$ is
chosen so that $|r_0|^2=1$ at some initial point $0\in\Z^2$.
Now we obtain a purely imaginary pair
$$
  \tilde\sigma^\pm
  = {\pm i\over|r\mp i|^2}\{
    (\sigma^++\sigma^-) \pm i\,(r\sigma^+-\bar r\sigma^-)\},
$$
by choosing $c^\pm=\pm i$ again.
Hence $\mp i\tilde\sigma^\pm$ define K\"onigs dual lifts%
\footnote{Note that reciprocal constant rescaling of
  $\sigma^\mp$ demands rescaling of $r$
  in \reqn{christoffel} by the same factor.
  In the case at hand, the purely imaginary function
  $\tilde r$ is turned into the real function $i\tilde r$.}
of a real $\Omega$-net and we conclude:

\proclaim\htag[thm.cplx]{\thmno Cor}.
Any $\Omega$-net can be spanned by pairs of isothermic sphere
congruences with complex conjugate K\"onigs dual lifts and,
conversely, any such pair gives rise to an $\Omega$-net.

Thus our \href[def.omega]{Def 3.1} encompasses also the
linear Weingarten nets with $\alpha\gamma-\beta^2>0$
by \href[thm.lwo]{Thm 2.8}.

The transformations of the enveloped isothermic sphere
congruences of an $\Omega$-net give rise to transformations
of the net.
As the isothermic transformation theory of discrete isothermic
nets hinges on the isothermic loop of flat connections of the
net which, in turn, depends on the cross ratio factorizing
function of the net, we shall start by getting our hands on
this function.

First observe that rewriting \reqn{christoffel} we obtain
(a lift of) the edge curvature sphere
$$
  \kappa_{ij}
  = r_ir_j\sigma^+_i - \sigma^-_i
  = r_ir_j\sigma^+_j - \sigma^-_j
  \in f_i\cap f_j
\eqno cslift$$
of the Legendre map $f$ on an edge $(ij)$, cf \reqn{cs}.
Using that $\kappa_{ij}\perp\sigma^\pm_i,\sigma^\pm_j$
we learn that
$$
    (\sigma^-_i\sigma^+_j)
  = (\sigma^+_i\sigma^-_j)
  = (\mu^\pm_i\mu^\pm_j)
  =: a_{ij}.
\eqno elab$$
Clearly, $a$ is an edge function, $a_{ij}=a_{ji}$.
Rearranging the Moutard equation for $\mu^+$ from \reqn{kmeq}
suitably and taking norm squares,
$$
    {\mu^+_k\over r_k-r_i} - {\mu^+_l\over r_l-r_j}
  = {\mu^+_i\over r_k-r_i} - {\mu^+_j\over r_l-r_j}
    \enspace\Rightarrow\enspace
  {a_{ij}-a_{kl}\over(r_k-r_i)(r_l-r_j)} = 0,
$$
we also learn that $a$ is an {\em edge-labelling\/}
in the sense of \href[ref.bosu08]{[2, Def 4.4]},
that is, it is constant across opposite edges of faces,
$a_{ij}=a_{kl}$,
cf \href[ref.bosu08]{[2, Thm 4.29]}.

To see that $a$ is indeed a cross ratio factorizing function,
cf \href[ref.imdg]{[14, \S5.7.2]}
or \href[ref.do06]{[12, Prop 10]},
first note that the vertices of a face of either isothermic
sphere congruence $s^\pm$ lie on a conic in a projective
plane since $s^\pm$ are Ribaucour sphere congruences in
the sense of \href[ref.bosu08]{[2, Def 3.27]}.
Fixing three points $s^\pm_i$, $s^\pm_j$ and $s^\pm_l$
of a face, the cross ratio
$q=[s^\pm_i,s^\pm_j,s^\pm_k,s^\pm_l]\in\R\cup\{\infty\}$
bijectively parametrizes the conic via
$$
  s^\pm_k
  = \langle \sigma^\pm_i + {1\over(\sigma^\pm_j\sigma^\pm_l)}\{
      (q-1)(\sigma^\pm_i\sigma^\pm_l)\sigma^\pm_j
    + ({1\over q}-1)(\sigma^\pm_i\sigma^\pm_j)\sigma^\pm_l\}
    \rangle,
\eqno diaconn$$
where $\sigma^\pm$ is any lift of $s^\pm$,
cf \href[ref.do06]{[12, (B.7)]}%
\footnote{Note that
  $
    q = [s^\pm_i,s^\pm_j,s^\pm_k,s^\pm_l]
    = cr(s^\pm_j,s^\pm_l,s^\pm_i,s^\pm_k)
  $ as
  our definition of the cross ratio differs from the classical
  one used in \href[ref.do06]{[12]} by the order of points.}
or \href[ref.bjrs08]{[6, Sect 2.1]}.
It is now straightforward to verify that
$q={a_{ij}\over a_{jk}}$ as, for Moutard
lifts $\mu^\pm$ of $s^\pm$ and taking inner
products with $\mu^\pm_j$ in \reqn{kmeq},
$$
  \mu^\pm_k
  = \mu^\pm_i
  - {a_{ij}-a_{jk}\over(\mu^\pm_j\mu^\pm_l)}\delta\mu^\pm_{jl}
  = \mu^\pm_i
  + {r^{\pm1}_k-r^{\pm1}_i\over r^{\pm1}_l-r^{\pm1}_j}
    \delta\mu^\pm_{jl}.
\eqno ameq$$

In summary:
\proclaim\htag[thm.crfact]{\thmno Lemma}.
The edge-labelling $a$ of \reqn{elab} factorizes the
cross ratios of faces of either isothermic sphere congruence,
$$
  [s^\pm_i,s^\pm_j,s^\pm_k,s^\pm_l]={a_{ij}\over a_{jk}}.
$$
In particular, the edge-labelling $a$ of \reqn{elab} is,
up to constant rescaling,
a well defined Lie geometric invariant of each
isothermic sphere congruence, $s^+$ and $s^-$.

We are now in a position to exploit the zero-curvature
representation of discrete isothermic nets.  We begin with a
rapid review of the formalism of metric connections on the discrete
vector bundle $\Z^2\times\R^{4,2}\to\Z^2$ cf
\href[ref.bjrs08]{[6, Def 2.4]}: a metric connection
$\Gamma$ on $\Z^2\times\R^{4,2}$ assigns to each
oriented edge $(ij)$, a linear isometry
$
\Gamma_{ij}:\{j\}\times\R^{4,2}\to\{i\}\times\R^{4,2}
$
such that $\Gamma_{ji}=\Gamma_{ij}^{-1}$, for all edges
$(ij)$.  In this context, a {\em gauge transformation\/} is
a map $i\mapsto T_i:\Z^2\to\SO(4,2)$ where we view $T_i$ as a linear
isometry of $\{i\}\times\R^{4,2}$.  Gauge transformations $T$
act on connections $\Gamma$ by
$$
(T\Gamma)_{ij}=T_i\Gamma_{ij}T_j^{-1}.
$$

A connection $\Gamma$ is {\em flat\/} if, on every
elementary quadrilateral $(ijkl)$, we have 
$$
  \Gamma_{ij}
  \Gamma_{jk}
  \Gamma_{kl}
  \Gamma_{li}
  = id
  \enspace\hbox{\rm or, equivalently,}\enspace
    \Gamma_{ij} \Gamma_{jk}
  = \Gamma_{il} \Gamma_{lk}.
$$
In this case, we can trivialise the connection, that is,
there is a gauge transformation $T$ with $T\Gamma=id$:
$$
\Gamma_{ij}=T_i^{-1}T_{j},
$$
for all edges $(ij)$.  Clearly, any gauge transform of a
flat connection is also flat.

With this understood, we are able to introduce the isothermic
loop of connections of an isothermic sphere congruence,

\proclaim\htag[def.isoconn]{\thmno Def}.
Let $s:\Z^2\to{\cal L}^4$ be an isothermic sphere congruence
with cross ratio factorizing edge-labelling $a$.
The {\em isothermic loop of connections\/}
$(\Gamma(t))_{t\in\R}$ of $s$ is a $1$-parameter family of
connections given by
$$
\Gamma_{ij}(t)x := \cases{
    (1-ta_{ij})\,x & if $x\in s_i$, \cr
    x & if $x\in(s_i\oplus s_j)^\perp$, \cr
    {1\over 1-ta_{ij}}\,x & if $x\in s_j$. \cr
  }
$$

Clearly $\Gamma_{ji}(t)\Gamma_{ij}(t)=id$ away from the
singularity $t={1\over a_{ij}}$, so that $\Gamma(t)$ defines
indeed a connection on the discrete vector bundle
$\Z^2\times\R^{4,2}$.
When $\sigma$ denotes any lift of the isothermic sphere
congruence $s$ then
$$
  \Gamma_{ij}(t)x
  = x + {ta_{ij}\over(\sigma_i\sigma_j)}\{
    {1\over 1-ta_{ij}}(x\sigma_i)\sigma_j
    - (x\sigma_j)\sigma_i\}.
$$
Note the structural similarity to \reqn{diaconn}
--- indeed parallel sections of $\Gamma(t)$ in the Lie quadric
yield Darboux transforms of $s$:
the corresponding condition on edges realizes a propagation
by cross ratio $ta$, thus yields a discrete version of
{\em Darboux's linear system\/},
cf \href[ref.bjrs08]{[6, Def 4.1]}.
For both, Darboux and Calapso transformations of an
isothermic sphere congruence, flatness of the connections
$\Gamma(t)$ is paramount.
Thus, returning to our context of an $\Omega$-net $f$
enveloped by a pair $s^\pm$ of isothermic sphere congruences,
we aim to convince ourselves that the connections $\Gamma^+(t)$
of $s^+$ are flat,
cf \href[ref.bjrs08]{[6, Lemma 2.5]}:

\proclaim\htag[thm.flatness]{\thmno Lemma}.
Given an $\Omega$-net $f=s^+\oplus s^-$ in terms of a pair of
isothermic sphere congruences $s^\pm$ that admit K\"onigs dual
lifts $\sigma^\pm$ the isothermic loop of connections of $s^+$
consists of flat connections.

Thus we wish to show that, on an elementary quadrilateral
$(ijkl)$ and as long as $t\neq{1\over a_{ij}},{1\over a_{jk}}$,
$$
    \Gamma^+_{ij}(t) \Gamma^+_{jk}(t)
  = \Gamma^+_{il}(t) \Gamma^+_{lk}(t).
$$
Having obtained the cross ratio factorizing property of the
edge-labelling $a$ in \href[thm.crfact]{Lemma 3.5} above,
the relevant part of the proof of
\href[ref.bjrs08]{[6, Lemma 2.5]}
applies, asserting correctness of the claim.

For autonomy we outline a simple algebraic proof here:
first observe that, clearly,
$$
  \Gamma^+_{ij}(t)\Gamma^+_{jk}(t)\,x
  = {1-ta_{jk}\over 1-ta_{ij}}\,x
    \enspace{\rm for}\enspace
  x \in s^+_j
$$
and a straightfoward computation, using the Moutard lift
$\mu^+=r\sigma^+$ of $s^+$ and \reqn{ameq}, shows that
$$
  \Gamma^+_{ij}(t)\Gamma^+_{jk}(t)\,x
  = {1-ta_{ij}\over 1-ta_{jk}}\,x
    \enspace{\rm for}\enspace
  x \in s^+_l.
$$
Since $\Gamma^+_{ij}(t)\Gamma^+_{jk}(t)\in\SO(\R^{4,2})$
acts trivially on $(s^+_i\oplus s^+_j\oplus s^+_l)^\perp$
it must act trivially on $(s^+_j\oplus s^+_l)^\perp$
and flatness of $\Gamma^+(t)$ follows by symmetry and
the fact that $a$ is an edge-labelling.%
\footnote{Note how, conversely, the limit $t\to\infty$
  yields the cross ratio factorizing nature
  of the edge-labelling $a$
  of \href[thm.crfact]{Lemma 3.5}.}

Instead of using a symmetry argument to deduce the flatness
of the connections $\Gamma^-(t)$ we employ a gauge theoretic
argument: we shall see that $\Gamma^-(t)$ and $\Gamma^+(t)$
are gauge equivalent, hence the flatness of $\Gamma^+(t)$
from \href[thm.flatness]{Lemma 3.7} implies flatness of
$\Gamma^-(t)$.

To this end let $g$ be any function on $\Z^2$ and consider
the following gauge transform of $\Gamma^-(t)$:
$$
  \Gamma^g_{ij}(t) :
  = (A^g\Gamma^-)_{ij}(t)
  = A^g_i(t)\Gamma^-_{ij}(t)A^g_j(-t),
    \enspace{\rm for}\enspace
  A^g(t):=1-tg\,(\sigma^+\wedge\sigma^-)
  = \exp(-tg\,\sigma^+\wedge\sigma^-),
\eqno gfam$$
where we identify ${\fk so}(4,1)\cong\bigwedge^2\R^{4,2}$
  via $(x\wedge y)z=(xz)y-(yz)x$.
Then any of the $\Gamma^g(t)$ is a metric connection on the
discrete vector bundle $\Z^2\times\R^{4,2}$.

Next note that the connections $\Gamma^g(t)$ have the
same shape as the connections $\Gamma^\pm(t)$ of the
isothermic loops of connections of $s^\pm$:
firstly, $\Gamma^g_{ij}(t)$ acts trivially on the
curvature sphere $f_i\cap f_j$;
secondly, using the lift \reqn{cslift} of the curvature sphere
$\langle\kappa_{ij}\rangle=f_i\cap f_j$, we learn that
$\Gamma^g_{ij}(t)$ has eigenspaces%
\footnote{At this point we see that the $\Gamma^g$ generally
  do not come from an isothermic sphere congruence:
  the condition that the eigenspaces $x_{ij}\subset f_i$
  coincide for all incident edges imposes a restriction
  on the function $g$.}
$$
  x_{ij} :
  = A^g_j({1\over a_{ij}})s^-_i
  = \langle\sigma^-_i + g_j\kappa_{ij}\rangle
    \enspace{\rm and}\enspace
  x_{ji} :
  = A^g_i({1\over a_{ij}})s^-_j
  = \langle\sigma^-_j + g_i\kappa_{ij}\rangle
$$
with eigenvalues $(1-ta_{ij})^{\pm1}$, respectively,
since
$
    (\sigma^+\wedge\sigma^-)_j\sigma^-_i
  = (\sigma^+\wedge\sigma^-)_i\sigma^-_j
  = -a_{ij}\kappa_{ij}
$;
and, finally, $\Gamma^g_{ij}(t)$ acts trivially on
$(f_i+f_j)^\perp$
--- hence
$$
  \Gamma^g_{ij}(t)x = \cases{
    (1-ta_{ij})\,x & if $x\in x_{ij}\subset f_i$, \cr
    x & if $x\in(x_{ij}\oplus x_{ji})^\perp$, \cr
    {1\over 1-ta_{ij}}\,x & if $x\in x_{ji}\subset f_j$. \cr
  }
$$

In particular, $x_{ij}=s^+_i$ and $x_{ji}=s^+_j$ for
$g\equiv 1$, showing that $\Gamma^+=\Gamma^1$ so that
flatness of $\Gamma^+(t)$ from \href[thm.flatness]{Lemma 3.7}
yields flatness of $\Gamma^-(t)=\Gamma^0(t)$,
hence of all connections $\Gamma^g(t)$:

\proclaim\htag[thm.gfamflat]{\thmno Cor}.
All connections $\Gamma^g(t)$ defined by \reqn{gfam} are flat.

Note that the gauge family of loops of connections $\Gamma^g$
also comprises the isothermic loops of the sphere congruences
$\tilde s^\pm=\langle\tilde\sigma^\pm\rangle$ of
\reqn{ambiguity}:
$\tilde s^\pm$ and $s^\pm$ share the same edge-labelling
$a$ and
$$
  x_{ij}
  = \langle\sigma^-_i+c^\pm r_i\sigma^+_i\rangle
  = \tilde s^\pm_i
    \enspace{\rm and}\enspace
  x_{ji}
  = \langle\sigma^-_j+c^\pm r_j\sigma^+_j\rangle
  = \tilde s^\pm_j
    \enspace{\rm for}\enspace
  g={c^\pm\over r+c^\pm}.
$$
In particular, when $s^\pm$ have a complex conjugate pair
of K\"onigs dual lifts, the gauge family contains real
connections.
A simple choice is given by $g\equiv{1\over 2}$:
using that $r_ir_j$ is unitary, $r_ir_j=e^{-2i\alpha}$,
we find that $\Gamma^g_{ij}(t)$ has real eigenspaces
$$
  x_{ij}
  = \langle e^{i\alpha}(r_ir_j\sigma^+_i+\sigma^-_i)\rangle
    \enspace{\rm and}\enspace
  x_{ji}
  = \langle e^{i\alpha}(r_ir_j\sigma^+_j+\sigma^-_j)\rangle.
$$

By their flatness all connections $\Gamma^g(t)$ can be
trivialized: there are (away from singularities)
maps
$$
  T^g(t):\Z^2\to\SO(\R^{4,2})
    \enspace\hbox{\rm so that}\enspace
  (T^g\Gamma^g)_{ij}(t)
  = T^g_i(t)\Gamma^g_{ij}(t)(T^g_j(t))^{-1}
  = id
\eqno triv$$
on every edge $(ij)$ of $\Z^2$.
And, as the connections are gauge equivalent via \reqn{gfam},
the gauge transformations $T^g(t)$ are, up to constants of
integration, related by
$
  T^-(t) = (T^gA^g)(t).
$
In particular, recall that the {\em Calapso
transformation\/} $T$ of an isothermic surface (in a
    quadric of any signature) is obtained by trivialising
    the loop of isothermic connections, see
    \href[ref.bjrs08]{[6, Def 2.7]} or \href[ref.bdpp11]{[5,
    Thm 4.15]} where this is
    treated in a rather general setting.  In the present
    case, the Calapso transformations $T^{\pm}$
of the pair $s^\pm$ of isothermic sphere congruences
are related by $T^-=T^+A^1$.
Thus, as $A^g_i(t)$ acts trivially on the contact element
$f_i$, any of the gauge transformations $T^g$ realizes the
{\em Calapso transforms\/} $T^\pm(t)s^\pm=T^g(t)s^\pm$
of both isothermic sphere congruences $s^\pm$ spanning
an $\Omega$-net.
This motivates the following definition,
cf \href[ref.bjr10]{[7, Sect 3]}
and \href[ref.bjrs08]{[6, Lemma 2.7]}:

\proclaim\htag[def.calapso]{\thmno Thm \& Def}.
Let $f=s^+\oplus s^-$ be an $\Omega$-net spanned by a pair of
isothermic sphere congruences $s^\pm$ admitting K\"onigs dual
lifts $\sigma^\pm$ and
let $T^g(t)$ be trivializing gauge transformations \reqn{triv}
of the connections $\Gamma^g(t)$ of \reqn{gfam}.
Then, the deformation
$$
  t \mapsto f(t) := T^g(t)f = T^\pm(t)f
$$
does not depend on the choice of gauge function $g$.
It will be called the {\em Calapso deformation\/} of the
$\Omega$-net $f$.
The {\em Calapso transforms\/} $f(t)$ of $f$ are $\Omega$-nets
with enveloped isothermic sphere congruences
$s^\pm(t)=T^g(t)s^\pm=T^\pm(t)s^\pm$
and K\"onigs dual lifts $\sigma^\pm(t)=T^g(t)\sigma^\pm$.

Only the last claim of the theorem
--- that $f(t)=(s^+\oplus s^-)(t)$ is an $\Omega$-net
    with $s^\pm(t)$ admitting K\"onigs dual lifts ---
requires further thought.
As the Calapso transform does not depend on the choice
of $g$ we may, without loss of generality,
assume $g\equiv 0$, that is, $T^g=T^-$.

Since $T^-(t):\Z^2\to\O(\R^{4,2})$ the Calapso transforms
$f(t)$ of $f$ take values in the space of contact elements.
To prove that $\sigma^\pm(t)$ are edge-parallel,
hence $f(t)$ satisfies the contact condition,
we employ the lift \reqn{cslift} of the curvature spheres:
then, on any edge $(ij)$,
$$
  r_ir_j d\sigma^+_{ij}(t)
  = T^-_i(t)\{\Gamma^-_{ij}(t)(\sigma^-_j+\kappa_{ij})
                            - (\sigma^-_i+\kappa_{ij})\}
  = T^-_i(t)\{\Gamma^-_{ij}(t)\sigma^-_j - \sigma^-_i\}
  = d\sigma^-_{ij}(t).
\eqno calchrist$$
Thus $\sigma^\pm(t)$ are edge-parallel and $f(t)$ is
a Legendre map with curvature spheres
$$
  \langle\kappa_{ij}(t)\rangle
  = \langle r_ir_j\sigma^+_i(t) - \sigma^-_i(t)\rangle
  = \langle r_ir_j\sigma^+_j(t) - \sigma^-_j(t)\rangle
  = f_i(t)\cap f_j(t).
$$
This also teaches us that $\mu^\pm(t):=r^{\pm1}\sigma^\pm(t)$
with $r(t)=r$ yields Moutard lifts of $s^\pm(t)$
since the Moutard equations \reqn{kmeq} were just
the integrability conditions of \reqn{christoffel}.
K\"onigs duality (away from umbilical faces)
of $\sigma^\pm(t)$, see \reqn{kmeq},
now also follows directly from \reqn{calchrist}:
$$\matrix{
  r_ir_k(r_l-r_j)\,\delta\sigma^+_{ik}(t)
  &=& r_ir_kr_l(d\sigma^+_{il}(t) + d\sigma^+_{lk}(t))
   -  r_ir_jr_k(d\sigma^+_{ij}(t) + d\sigma^+_{jk}(t)) \hfill\cr
  &=& r_k(d\sigma^-_{ji}(t) + d\sigma^-_{il}(t))
   -  r_i(d\sigma^-_{jk}(t) + d\sigma^-_{kl}(t)) \hfill\cr
  &=& (r_k-r_i)\,\delta\sigma^-_{jl}(t). \hfill\cr
}$$
This concludes the proof of \href[def.calapso]{Thm 3.9}.

Note that, in contrast to the function $r$ relating
K\"onigs dual and Moutard lifts that is invariant
under the Calapso deformation,
the edge-labelling $a$ changes,
cf \href[ref.imdg]{[14, \S5.7.16]}
or \href[ref.bjrs08]{[6, Lemma 2.7]}:
$$
  a_{ij}(t)
  = (\mu^\pm_i(t)\mu^\pm_j(t))
  = (\mu^\pm_i\Gamma^\pm_{ij}\mu^\pm_j)
  = {1\over 1-ta_{ij}}(\mu^\pm_i\mu^\pm_j)
  = {a_{ij}\over 1-ta_{ij}}.
$$

\section Lawson transformation of linear Weingarten surfaces.

In this section, we shall see how the Calapso deformation
of \href[def.calapso]{Thm 3.9} for $\Omega$-nets descends
to a ``Lawson transformation'' for linear Weingarten nets.
To this end we need to investigate the effect of the
deformation on the two linear sphere complexes that
come with a linear Weingarten net.
Thus let $f=s^+\oplus s^-$ be an $\Omega$-net so that $s^\pm$
have K\"onigs dual lifts $\sigma^\pm$ and take values in linear
sphere complexes ${\fk k}^\pm$, where we assume the above
relative normalizations $({\fk k}^\pm\sigma^\mp)=-1$.

For symmetry we base our analysis on the ``middle connection''
$\Gamma^g$ with $g\equiv{1\over 2}$.
Recall that $\Gamma^g$ is then real in the complex conjugate
case and, consequently, so is $T^g$ when fixing $T^g$ to
be real at an initial point,
say $T^g_0=id$ at some point $0\in\Z^2$.
Now
$$
  {\fk k}^\pm(t) :
  = T^g(t)\{{\fk k}^\pm + {t\over 2}\,\sigma^\pm \}
  = T^\pm(t){\fk k}^\pm
  \equiv const
\eqno kpmt$$
for any fixed $t$ since
${\fk k}^\pm\perp s^\pm_i\oplus s^\pm_j$
so that $\Gamma^\pm_{ij}(t){\fk k}^\pm={\fk k}^\pm$
on any edge $(ij)$.
Moreover,
$$
  ({\fk k}^\pm(t)s^\pm(t)) \equiv 0
    \enspace{\rm for}\enspace
  s^\pm(t) = T^g(t)s^\pm
$$
since $T^g:\Z^2\to\O(\R^{4,2})$.
Thus ${\fk k}^\pm(t)$ define linear sphere complexes that
the isothermic sphere congruences $s^\pm(t)$ spanning the
Calapso transform $f(t)=(s^+\oplus s^-)(t)$ of $f$ take
values in.

Note that, with the K\"onigs dual lifts
$\sigma^\pm(t)=T^g(t)\sigma^\pm$ of $s^\pm(t)$,
see \href[def.calapso]{Thm 3.9}, the deformation preserves
the relative scaling $({\fk k}^\pm(t)\sigma^\mp(t))\equiv -1$
and, in the complex conjugate case, $\sigma^\pm(t)$ as well as
${\fk k}^\pm(t)$ are complex conjugate again.

Consequently, we have proved that the Calapso deformation
yields a transformation for linear Weingarten nets:
given a linear Weingarten net,
its Legendre lift is an $\Omega$-net
(by \href[thm.lwo]{Thm 2.8})
admitting (see \href[def.calapso]{Thm 3.9})
  Calapso deformation into a new $\Omega$-net,
which has the characteristics of the Legendre lift of a
  linear Weingarten net (cf \href[thm.lwo]{Thm 2.8}).
The only potential issue is that the resulting $\Omega$-net
may not admit an appropriate space form projection ---
if the plane spanned by ${\fk k}^+(t)$ and ${\fk k}^-(t)$
becomes null then it does not contain a point sphere
complex ${\fk p}(t)$.
Computing inner products we find
$$
  ({\fk k}^\pm(t){\fk k}^\pm(t)) = ({\fk k^\pm k^\pm})
    \enspace{\rm and}\enspace
  ({\fk k}^+(t){\fk k}^-(t)) = ({\fk k^+k^-}) - t,
\eqno linv$$
showing that this issue does generically not occur:
it only occurs at a single value of $t$ when ${\fk k}^\pm$
define a (possibly complex conjugate) pair of spheres.
Below we shall discuss these cases.

\proclaim\htag[def.lawson]{\thmno Thm \& Def}.
Let $f$ be the Legendre lift of a linear Weingarten net
$({\fk f,t}):\Z^2\to{\fk Q}^3\times{\fk P}^3$.
The Calapso transforms $f(t)$ of $f$ are generically%
\footnote{That is, as long as the sphere complexes
  ${\fk k}^\pm(t)$ from \reqn{kpmt} do not span
  a contact element.}
the Legendre lifts of suitable linear Weingarten nets
$({\fk f,t})(t):\Z^2\to{\fk Q}^3(t)\times{\fk P}^3(t)$.
These will be called {\em Lawson transforms\/}
of $({\fk f,t})$.

To justify the terminology we consider constant mean curvature
nets, cf \href[ref.bjrs08]{[6, Sect 5]}.

\htag[expl.cmc]{{\bf\thmno Example}}.
Let $({\fk f,t}):\Z^2\to{\fk Q}^3\times{\fk P}^3$ be a constant
mean curvature net, that is, a linear Weingarten net with
$$
  0 = \alpha\,K + 2\beta H + \gamma,
    \enspace{\rm where}\enspace
  \alpha = 0.
\eqno lwcmc$$
In the discussion leading up to \href[thm.lwo]{Thm 2.8}
we already saw, cf \href[ref.bjl11]{[4, Lemma 4.1]},
that 
$$
  \sigma^- = {\fk f}
    \enspace{\rm and}\enspace
  \sigma^+ = {\fk t} + H\,{\fk f}
$$
yield a K\"onigs dual pair of isothermic sphere congruences
$s^\pm$ that take values in linear sphere complexes
$$
  {\fk k}^- = {\fk p}
    \enspace{\rm and}\enspace
  {\fk k}^+ = {\fk q} - H\,{\fk p}.
$$

To recover the Lawson correspondence
of \href[ref.bjrs08]{[6, Def 5.2]}
we follow the arguments that proved
\href[def.lawson]{Thm 4.1} but
now base our analysis on $T^-$ instead of $T^{1\over 2}$.
Since $\Gamma^-_{ij}(t){\fk k}^-={\fk k}^-$ for all $t$ and
all edges $(ij)$ we may assume that
$T^-(t){\fk p=p}$ for all $t$,
that is, $T^-(t)$ is a M\"obius geometric Calapso transformation
of the discrete isothermic net ${\fk f}:\Z^2\to{\fk Q}^3$
when $({\fk pp})=-1$,
see \href[ref.imdg]{[14, \S5.7.16]}.
Hence, for all $t$,
$$
  {\fk k}^-(t) = {\fk k}^- = {\fk p}
$$
provides a canonical point sphere complex for space form
projection of the Calapso transform $f(t)$ of the
Legendre lift $f$ of $({\fk f,t})$.
Further
$$
  {\fk k}^+(t)
  = T^-(t)\{{\fk k}^+ + t\sigma^+\}
$$
is obtained as the image of (the Lie geometric lift of)
the {\em linear conserved quantity\/}
of \href[ref.bjrs08]{[6, Def 5.1]}.
Consequently,
$$
  {\fk q}(t) :
  = {\fk k}^+(t) - {({\fk k}^+(t){\fk p})\over({\fk pp})}{\fk p}
  = T^-(t)\{{\fk q} + t(\sigma^++{\fk{p\over(pp)}})\}
$$
yields a canonical space form vector for the
space form projection.

Hence, as $\sigma^\pm(t)=T^-(t)\sigma^\pm$,
the corresponding space form projection is given by
$$
  {\fk f}(t)
  = T^-(t)\,{\fk f}:\Z^2\to{\fk Q}^3(t)
    \enspace{\rm and}\enspace
  {\fk t}(t)
  = T^-(t)\{{\fk t}-{t\over{\fk(pp)}}{\fk f}\}:
    \Z^2\to{\fk P}^3(t).
$$
Now $\sigma^+(t)={\fk t}(t)+(H+{t\over{\fk(pp)}})\,{\fk f}(t)$,
showing that the linear Weingarten net $({\fk f,t})(t)$ has
constant mean curvature
$$
  H(t)=H+{t\over{\fk(pp)}}.
$$

Note that we also recover the {\em Lawson invariant\/},
relating the mean curvature of the constant mean curvature net
and its ambient constant sectional curvature, cf \reqn{linv}:
$$
  ({\fk pp})\,H^2(t) + {\fk(q(t)q(t))}
  = {\fk(k^+(t)k^+(t))}
  = {\fk(k^+k^+)}
  = ({\fk pp})\,H^2 + {\fk(qq)}.
$$

Thus our Lawson transformation of \href[def.lawson]{Def 4.1}
does indeed generalize the Lawson correspondence
of \href[ref.bjrs08]{[6, Sect 5]}
for discrete constant mean curvature nets.

Coming back to the genericity issue of the
Lawson transformation from \href[def.lawson]{Thm 4.1}
we consider the discrete analogue
of (intrinsically) flat surfaces in hyperbolic space,
cf \href[ref.hrsy09]{[15, Sect 4.3]}
and \href[ref.bjr10]{[7]}.

\htag[expl.ff]{{\bf\thmno Example}}.
Thus let $({\fk f,t}):\Z^2\to{\fk Q}^3\times{\fk P}^3$,
where $({\fk pp})=-1$ and $({\fk qq})=+1$ so that ${\fk Q}^3$
becomes hyperbolic space, satisfy $A({\fk t,t})=A({\fk f,f})$,
cf \href[ref.hrsy09]{[15, Lemma 6.5]},
that is, $({\fk f,t})$ is linear Weingarten with
$$
  0 = \alpha\,K + 2\beta\,H + \gamma,
    \enspace{\rm where}\enspace
  \alpha + \gamma = \beta = 0.
\eqno lwff$$
Now $A({\fk f+t,f-t})=0$ so that a K\"onigs dual pair
of enveloped isothermic sphere congruences is given by
$$
  \sigma^\pm := {\fk f} \pm {\fk t}.
$$
As ${\fk f}$ and ${\fk t}$ take values
in ${\fk Q}^3$ and ${\fk P}^3$,
respectively, $s^\pm$ take values
in linear sphere complexes
$$
  {\fk k}^\pm := {1\over 2}({\fk q\mp p}).
$$
These define two oriented spheres since $({\fk k^\pm k^\pm})=0$:
they are the two orientations of the infinity boundary
of the ambient hyperbolic space,
cf \href[ref.bjr10]{[7, Sect 2]},
and the fact that $\sigma^\pm$ take values in the
sphere complexes
defined by ${\fk k}^\pm$ teaches us that $\sigma^\pm$ touch the
infinity boundary (with opposite orientations),
that is, $\sigma^\pm$ are {\em horosphere congruences\/}.

Before proceeding to the Lawson transformation,
note how the parallel nets of $({\fk f,t})$ are
obtained by simultaneous reciprocal rescaling
of ${\fk k}^\pm$ and $\sigma^\pm$,
cf \reqn{cbas} and \href[thm.lwparallel]{Cor 2.9}:
with
$$
  \tilde{\fk k}^\pm = e^{\pm\rho}{\fk k}^\pm
    \enspace{\rm and}\enspace
  \tilde\sigma^\pm  = e^{\pm\rho}\sigma^\pm
$$
the relative scalings are preserved and a new choice
of point sphere complex and space form vector
$$
  \left.\matrix{
    \tilde{\fk p} :
    = \tilde{\fk k}^- - \tilde{\fk k}^+
    = {\fk p\cosh\rho - q\sinh\rho} \hfill\cr
    \tilde{\fk q} :
    = \tilde{\fk k}^- + \tilde{\fk k}^+
    = {\fk q\cosh\rho - p\sinh\rho} \hfill\cr
  }\right\}
    \enspace{\rm yields}\enspace
  \left\{\matrix{
    \tilde{\fk f}
    = {\fk f\cosh\rho + t\sinh\rho}, \hfill\cr
    \tilde{\fk t}
    = {\fk f\sinh\rho + t\cosh\rho}. \hfill\cr
  }\right.
$$
The linear Weingarten condition \reqn{lwff} is preserved
by this change of space form projection by parallel
transformation.
Thus the parallel nets of a flat net in hyperbolic space
are flat.
For the analysis of the Lawson transformation below
we shall disregard this freedom.

To honour the symmetry of the situation we again use the
Calapso transformations $T^g(t)$ of the ``middle connection''
$\Gamma^g(t)$ with $g\equiv{1\over 2}$, as above
and in \href[ref.bjr10]{[7, Sect 3]}.
Thus, cf \reqn{kpmt}, we obtain
$$
  \sigma^\pm(t)
  = T^g(t)\sigma^\pm
  = T^g(t)\{{\fk f\pm t}\}
    \enspace{\rm and}\enspace
  {\fk k}^\pm(t)
  = T^g(t)\{{\fk k}^\pm + {t\over 2}\sigma^\pm\}
  = T^g(t)\{
    {{\fk q}+t{\fk f}\over 2} \mp {{\fk p}-t{\fk t}\over 2}
    \}
$$
as the K\"onigs dual pair spanning the Calapso transform
$f(t)$ of the original $\Omega$-net and the sphere complexes
the enveloped isothermic sphere congruences take values in.
By \reqn{linv} the new sphere complexes still define two
oriented spheres, $({\fk k^\pm k^\pm})=0$, which can be
interpreted as the two orientations of the infinity boundary
of a hyperbolic space as long as they do not touch,
$$
  ({\fk k^+(t)k^-(t)})
  = ({\fk k^+k^-}) - t
  = {1\over 2}(1-2t)
  \neq 0,
$$
that is, as long as the genericity condition
of \href[def.lawson]{Thm 4.1} is satisfied.%
\footnote{When $t={1\over 2}$ the plane
  $\langle{\fk k^+,k^-}\rangle$ is isotropic
  hence does not contain a point sphere complex.}
Then a choice of
$$
  {\fk p}(t) :
  = {\fk k}^-(t) - {\fk k}^+(t)
  = T^g(t)({\fk p}-t\,{\fk t})
    \enspace{\rm and}\enspace
  {\fk q}(t) :
  = {\fk k}^-(t) + {\fk k}^+(t)
  = T^g(t)({\fk q}+t\,{\fk f})
$$
for the point sphere complex and space form vector yields a
projection to a hyperbolic space and a de Sitter space as its
(unit) tangent bundle, or vice versa, depending on the sign
of $1-2t$:
$$
  ({\fk p}(t){\fk p}(t)) = -(1-2t)
    \enspace{\rm and}\enspace
  ({\fk q}(t){\fk q}(t)) = +(1-2t).
$$
The corresponding space form projection
$$
  ({\fk f,t})(t)
  = (T^g(t){\fk f},T^g(t){\fk t}):
    \Z^2\to{\fk Q}^3(t)\times{\fk P}^3(t)
$$
is a linear Weingarten net satisfying the same linear Weingarten
condition \reqn{lwff} since, as before,
$$
  W(t) :
  = 2\,{\fk k}^+(t)\odot{\fk k}^-(t)
  = {1\over 2}\{{\fk q}(t)\odot{\fk q}(t)
  - {\fk p}(t)\odot{\fk p}(t)\}.
$$
In particular, as long as $t<{1\over 2}$, the Lawson transforms
$({\fk f,t})(t)$ of $({\fk f,t})$ remain discrete analogues
of (intrinsically) flat surfaces in hyperbolic space.%
\footnote{Rescaling ${\fk p}(t)$ and ${\fk q}(t)$ by
  ${1\over\sqrt{1-2t}}$ yields the standard model of
  hyperbolic space with de Sitter space as its unit
  tangent bundle.}
Beyond the singularity of the Lawson transformation,
when $t>{1\over 2}$,
we obtain linear Weingarten nets $({\fk f,t})(t)$
with constant (extrinsic) Gauss curvature $K=1$
in de Sitter space.%
\footnote{Note that, swapping the roles of ${\fk f}(t)$
  and ${\fk t}(t)$ in this case, yields linear Weingarten
  nets in hyperbolic space again.}

Thus we obtained a case, where the genericity issue for the
Lawson transformation does occur and,
in particular,
we have seen how the two conserved quantities ${\fk k}^\pm(t)$
in this case become spheres
--- the two orientations of the infinity spheres
    of the ambient hyperbolic geometries ---
as predicted from \reqn{linv}.
Below we shall give a more exhaustive discussion
of the genericity phenomenon.

Generalizing \href[expl.ff]{Example 4.3} we next discuss
discrete nets of (arbitrary) constant Gauss curvature in
(possibly Lorentzian) space forms:

\htag[expl.cgc]{{\bf\thmno Example}}.
Let $({\fk f,t}):\Z^2\to{\fk Q}^3\times{\fk P}^3$,
where $\varepsilon:=-({\fk pp})=\pm1$ and $\kappa=-({\fk qq})$,
be a space form projection into a quadric ${\fk Q}^3$
of constant curvature $\kappa$ satisfying
$$
  0 = \alpha\,K + 2\beta\,H + \gamma
    \enspace{\rm with}\enspace
  \beta = 0.
\eqno lwcgc$$
Clearly,%
\footnote{Otherwise $\gamma=0$ as well and the linear
  Weingarten condition would become trivial.}
$\alpha\neq0$ so that we may, without loss of generality,
assume that $\alpha=1$;
excluding tubular linear Weingarten surfaces we have
$K=-\gamma\neq0$.
Now
$$
  2W = 4\,{\fk k^+\odot k^-}
  = {\fk p\odot p} - {1\over K}\,{\fk q\odot q}
    \enspace{\rm with}\enspace
  {\fk k}^\pm :
  = {1\over 2}\{{\fk p} \mp {1\over\sqrt{K}}\,{\fk q}\},
$$
showing that ${\fk k}^\pm$ and the corresponding
K\"onigs dual lifts
$\sigma^\pm={\fk t}\pm\sqrt{K}\,{\fk f}$
of enveloped isothermic sphere congruences $s^\pm$
become complex conjugate when $K<0$.

To investigate the Lawson transformation, we use the
``middle connection'' $\Gamma^g$ with $g\equiv{1\over 2}$
again, as in the case of flat fronts in hyperbolic space:
recall that this connection is real in the complex
conjugate case so that the Calapso transformations
$T^g(t)$ can be assumed to be real as well.
Now, cf \reqn{kpmt},
$$
  \sigma^\pm(t)
  = T^g(t)\{{\fk t}\pm\sqrt{K}\,{\fk f}\}
    \enspace{\rm and}\enspace
  {\fk k}^\pm(t)
  = T^g(t)\{
    {{\fk p}+t{\fk t}\over 2}
    \mp {{\fk q}-tK{\fk f}\over 2\sqrt{K}}\}
$$
yields the Calapso transform $f(t)$ of the original
$\Omega$-net $f=({\fk f,t})$ and the two sphere complexes
that its pair of enveloped isothermic sphere congruences
$s^\pm(t)=\langle\sigma^\pm(t)\rangle$ take values in.
Thus, as long as ${\fk p}+t{\fk t}$ does not become isotropic,
$\varepsilon+2t\neq0$, we may choose
$$
  {\fk p}(t) :
  = {1\over\sqrt{|\varepsilon+2t|}}
    T^g(t)\{{\fk p} + t{\fk t}\}
    \enspace{\rm and}\enspace
  {\fk q}(t) :
  = T^g(t)\{{\fk q} - tK{\fk f}\}
$$
as the new point sphere complex and space form vector for
the space form projection $({\fk f,t})(t)$ of $f(t)$,
so that
$$
  {\fk f}(t) = T^g(t){\fk f}
    \enspace{\rm and}\enspace
  {\fk t}(t) = \sqrt{|\varepsilon+2t|}\,T^g(t){\fk t}.
$$
Note that $\varepsilon(t)=-({\fk p}(t){\fk p}(t))$ changes
sign at $t=-{\varepsilon\over 2}$, hence the ambient
geometry of $({\fk f,t})(t)$ changes signature,
as in the case of flat nets in hyperbolic space.

Now
$$
  2W(t)
  = 4{\fk k}^+(t)\odot{\fk k}^-(t)
  = |\varepsilon+2t|\,{\fk p}(t)\odot{\fk p}(t)
  - {1\over K}\,{\fk q}(t)\odot{\fk q}(t)
$$
encodes the linear Weingarten condition for $({\fk f,t})(t)$:
the new (constant) Gauss and ambient curvatures become
$$
  K(t) = K\,|\varepsilon+2t|
    \enspace{\rm and}\enspace
  \kappa(t) = -({\fk q}(t){\fk q}(t)) = \kappa - 2tK.
$$
As a consequence the ``intrinsic Gauss curvatures''
$$
  \gK(t) = \varepsilon(t)K(t) + \kappa(t) \equiv \gK
$$
of the nets remain unchanged by the Lawson transformation.
Of course, this fact depends on the chosen normalization
of the space form vectors ${\fk q}(t)$:
a rescaling of ${\fk q}(t)$ results in a rescaling of both
the extrinsic and intrinsic Gauss curvatures by the square
of the factor
---
in particular,
if ${1\over\sqrt{|\varepsilon+2t|}}{\fk q}(t)$
had been chosen for the space form projections instead,
the extrinsic Gauss curvatures would remain unchanged
while the intrinsic Gauss curvatures would get scaled
in the family.

Similar thoughts show that the apparent problem
of the space form projection when $t=-{\varepsilon\over 2}$
is easily resolved as long as the intrinsic Gauss curvature
of the original net did not vanish:
as long as neither ${\fk p}(t)$ nor ${\fk q}(t)$ become
isotropic their roles can be interchanged after suitable
rescalings.
In particular,
a common rescaling by ${1\over\sqrt{|\kappa(t)|}}$ and
re-interpretation of ${1\over\sqrt{|\kappa(t)|}}{\fk q}(t)$
as the point sphere complex results in a linear
Weingarten net $\sqrt{|\kappa(t)|}({\fk t,f})(t)$
of constant Gauss curvature ${1\over K(t)}$.
Hence the Lawson transformation is well defined as long as
${\fk p}(t)$ and ${\fk q}(t)$ do not simultaneously become
isotropic, that is, as long as
$$
  \kappa(-{\varepsilon\over 2}) = \gK \neq 0.
$$

As a third special class of linear Weingarten surfaces,
obtained by the vanishing of the third coefficient in
the linear Weingarten condition,
we discuss nets of constant harmonic mean curvature as
a discrete analogue of the surfaces with constant average
of their curvature radii:

\htag[expl.chmc]{{\bf\thmno Example}}.
Fix a space form projection
$({\fk f,t}):\Z^2\to{\fk Q}^3\times{\fk P}^3$,
where $-({\fk pp})=:\varepsilon=\pm1$ and the ambient
curvature is given by $-({\fk qq})=\kappa$ as before,
and suppose that
$$
  0 = \alpha\,K + 2\beta\,H + \gamma,
    \enspace{\rm where}\enspace
  \gamma = 0.
\eqno lwchmc$$
As we exclude tubular linear Weingarten nets we must have
$\beta\neq0$.
Hence, without loss of generality,
$\beta=1$ and $\alpha=-2{H\over K}$ is given by
the constant harmonic mean curvature of the net.
The enveloped isothermic sphere congruences of the constant
harmonic mean curvature net then turn out to be its tangent
plane congruence and its ``middle sphere congruence'',
$$
  \sigma^- = {\fk t}
    \enspace{\rm and}\enspace
  \sigma^+ = {\fk f} + {H\over K}\,{\fk t},
$$
which take values in the linear sphere complexes
$$
  {\fk k}^- = {\fk q}
    \enspace{\rm and}\enspace
  {\fk k}^+ = {\fk p} - {H\over K}\,{\fk q}.
$$

Note the similarity to the constant mean curvature nets
of \href[expl.cmc]{Example 4.2}.
Writing \reqn{lwchmc} in the more symmetric form \reqn{malw},
$$
  0 = \alpha\,A({\fk t,t})
    - 2\beta\,A({\fk t,f})
    + \gamma\,A({\fk f,f}),
$$
makes this similarity more tangible:
as long as $\kappa\neq0$ the aforementioned ``duality'' for
linear Weingarten nets relates constant harmonic mean curvature
nets and constant mean curvature nets.
In particular, a common rescaling of the point sphere complex
and the space form vector,
$$
  (\tilde{\fk q},\tilde{\fk p})
  = {1\over\sqrt{|\kappa|}}({\fk p,q}),
    \enspace{\rm yields}\enspace
  (\tilde{\fk f},\tilde{\fk t})
  = \sqrt{|\kappa|}\,({\fk t,f})
$$
as a constant mean curvature $\tilde H={H\over K}$ net
in a quadric of constant curvature
$\tilde\kappa={\varepsilon\over|\kappa|}$,
whose signature is given by
$\tilde\varepsilon={\kappa\over|\kappa|}$.
The Lawson transformation of \href[expl.cmc]{Example 4.2}
then yields a transformation into constant mean curvature
nets with%
\footnote{We use the edge-labelling $a$ of the original
  K\"onigs dual pair $(\sigma^+,\sigma^-)$ for the Calapso
  transformation here, which results in a rescaling of the
  spectral parameter of \href[expl.cmc]{Example 4.2}.}
$$
  \tilde H(t)
  = \tilde H - \tilde\varepsilon\,{t\over|\kappa|}
    \enspace{\rm and}\enspace
  \tilde\kappa(t)
  = \tilde\kappa
  + 2\tilde H\,{t\over|\kappa|}
  - \tilde\varepsilon\,({t\over|\kappa|})^2.
$$
As long as $\tilde\kappa(t)\neq0$ the same ``duality'' can
then be used to obtain constant harmonic mean curvature nets
$({\fk f}(t),{\fk t}(t))
 = \sqrt{|\tilde\kappa(t)|}(\tilde{\fk t}(t),\tilde{\fk f}(t))$
as Lawson transforms of the original net $({\fk f,t})$.

Aiming to obtain the Lawson transformation for nets of
constant harmonic mean curvature directly, we recover
the same regularity issues as outlined above.
Motivated by the observation that the characterizing
feature of constant harmonic mean curvature nets in our
setup is its tangent plane congruence ${\fk t}$ being one
of the enveloped isothermic sphere congruences,
we base our analysis on $T^-$, where constants of
integration are adjusted so that $T^-(t){\fk q}={\fk q}$.
Then we aim to obtain
$$
  {\fk k}^+(t)
  = T^-(t)\{{\fk k}^+ + t\sigma^+\}
  \parallel {\fk p}(t) - {H(t)\over K(t)}\,{\fk q}
$$
with a normalized point sphere complex ${\fk p}(t)\perp{\fk q}$
in order to recover the linear Weingarten condition of
a constant harmonic mean curvature net.
Orthogonalization then requires $\kappa\neq0$ and
normalization requires also $\tilde\kappa(t)\neq0$:
when both are satisfied%
\footnote{Note that $\tilde\kappa(t)\neq0$ encodes the
  non-degeneracy of the induced metric of the plane
  spanned by ${\fk k}^\pm(t)$.}
$$
  {\fk p}(t)
  = {1\over\sqrt{|\kappa\tilde\kappa(t)|}}\,T^-(t)
    \{ {\fk p} + t\,(\sigma^+ - {1\over\kappa}{\fk q}) \}
    \enspace{\rm yields}\enspace
  {1\over\sqrt{|\kappa\tilde\kappa(t)|}}{\fk k}^+(t)
  = {\fk p}(t)
    - {1\over\sqrt{|\kappa\tilde\kappa(t)|}}
      ({H\over K}-{1\over\kappa}\,t)\,{\fk q}.
$$

Starting from a constant harmonic mean curvature net in a
flat ambient geometry the Lawson transformation is still
well defined but the above approach,
following \href[expl.cmc]{Example 4.2},
does not yield a suitable space form projection.
However, exploiting the fact that constant harmonic
mean curvature nets in flat ambient geometries arise
as parallel nets of minimal nets, leads to an alternative
approach:
choosing
$$
  \tilde{\fk p}(t) := {\fk k}^+(t)
    \enspace{\rm and}\enspace
  \tilde{\fk q}(t) :
  = {\fk k}^-(t) - \varepsilon t\,{\fk k}^+(t)
$$
for a space form projection%
\footnote{Recall \reqn{linv}: when $\kappa=0$ we have
  $({\fk k}^+(t){\fk k}^+(t))=-\varepsilon$ and
  $({\fk k}^+(t){\fk k}^-(t))=-t$.}
we obtain a constant mean curvature net
$$
  (\tilde{\fk f}(t),\tilde{\fk t}(t))
  = (\sigma^+(t),\sigma^-(t)+\varepsilon t\,\sigma^+(t))
$$
as, clearly,
$
  A(\tilde{\fk f}(t),\tilde{\fk t}(t))
  = \varepsilon t\,A(\tilde{\fk f}(t),\tilde{\fk f}(t)).
$
Note that, in particular, this choice of projection yields
a minimal net $(\tilde{\fk f}(0),\tilde{\fk t}(0))$ for $t=0$.
When $t\neq0$ the ambient curvature of the constant mean
curvature net $(\tilde{\fk f}(t),\tilde{\fk t}(t))$ does
not vanish,
$\tilde\kappa(t)=-\varepsilon t^2\neq0$,
so that a choice
$$
  {\fk p}(t) := {1\over t}\tilde{\fk q}(t)
    \enspace{\rm and}\enspace
  {\fk q}(t) := \tilde{\fk p}(t)
$$
now yields a net
$({\fk f}(t),{\fk t}(t))=(\tilde{\fk t}(t),t\tilde{\fk f}(t))$
of constant harmonic mean curvature
${H(t)\over K(t)}=-\varepsilon$
in a (Lorentzian if $\varepsilon=1$) quadric
of constant curvature $\kappa(t)=\varepsilon$.

Indeed,
as our discussion of the Lawson transformation for constant
  harmonic mean curvature nets hints at,
the Lawson transformation is a transformation for parallel
families of linear Weingarten nets rather than for individual
nets: it involves a choice of space form projection that,
essentially, is a choice of a net in a parallel family.
We shall see that every parallel family of discrete linear
Weingarten nets contains nets of at least one of the
particular types discussed in the preceding examples.

\htag[expl.parallel]{{\bf\thmno Parallel families}}.
In \href[thm.lwparallel]{Cor 2.9} we already saw that parallel
nets of a linear Weingarten net are linear Weingarten,
cf \href[ref.imdg]{[14, Sect 2.7]}.
Using the same setup as above,
let $({\fk f,t})$ denote a linear Weingarten net in a space
form given by a point sphere complex ${\fk p}$ and a space
form vector ${\fk q\perp p}$.
A change of basis
$$
  (\tilde{\fk q},\tilde{\fk p}) = ({\fk q,p})B,
$$
where $B\in\Gl(2)$ is chosen to the preserve inner products
of the basis, yields a parallel linear Weingarten net
$(\tilde{\fk f},\tilde{\fk t})$.
The coefficients of the linear Weingarten relations \reqn{lw}
of $({\fk f,t})$ and $(\tilde{\fk f},\tilde{\fk t})$ are then
related by%
\footnote{This follows easily by interpreting the linear
  Weingarten condition \reqn{lw} as an orthogonality condition
  with respect to the inner product on symmetric
  $2\times2$-matrices given by the determinant
  as its quadratic form.}
$$
  \left(
    {\tilde\alpha\atop\tilde\beta}\,
    {\tilde\beta\atop\tilde\gamma}
  \right)
  = B
  \left(
    {\alpha\atop\beta}\,{\beta\atop\gamma}
  \right)
  B^t.
$$
As the shape of the basis transformations $B$ depends on the
signature of the plane $\langle{\fk q,p}\rangle$ we discuss
the cases that occur in turn,
cf \href[ref.pate88]{[17, Sect 3.4]}
or \href[ref.jtz97]{[13, Sect II.5]}.

\item{(1)} In the definite case
  we assume $({\fk q,p})$ to be an orthonormal basis
  so that the parallel family of nets is parametrized by
  $$
    B = \left(
      {\cos\vartheta\atop\sin\vartheta}\,
      {-\sin\vartheta\atop\phantom{-}\cos\vartheta}
      \right).
  $$
  Writing
  ${\alpha+\gamma\over 2}=\mu$,
  ${\alpha-\gamma\over 2}=\varrho\cos2\omega$ and
  $\beta=\varrho\sin2\omega$,
  the coefficients of the linear Weingarten condition \reqn{lw}
  of the parallel nets $({\fk f,t})(\vartheta)$ become
  $$
    \alpha(\vartheta) = \mu+\varrho\cos2(\vartheta+\omega),
      \enspace
    \beta(\vartheta)  = \varrho\sin2(\vartheta+\omega)
      \enspace{\rm and}\enspace
    \gamma(\vartheta) = \mu-\varrho\cos2(\vartheta+\omega).
  $$
\item{} Thus, if $\varrho\neq0$,
  then $\sin2(\vartheta+\omega)=0$ yields two pairs
  of antipodal constant Gauss curvature nets and,
  if also
  $$
    \mu^2 - \varrho^2 = \alpha\gamma-\beta^2 < 0,
  $$
  then $\cos2(\vartheta+\omega)=\mp{\mu\over\varrho}$
  yields two pairs of antipodal constant mean curvature nets
  and of constant harmonic mean curvature nets, respectively.
  When $\mu=0$ these coincide and yield minimal nets.
  Note the symmetric spacing of the twelve or eight nets
  that appear in this case.
\item{} If, on the other hand, $\varrho=0$
  then all $({\fk f,t})(\vartheta)$ are intrinsically flat,
  cf \href[expl.cgc]{Example 4.4}.

\item{(2)} In the degenerate case
  we assume that $({\fk qq})=0$ and $({\fk pp})=\pm1$
  so that the parallel family is parametrized by
  $$
    B = \left(
      {1\atop0}\,
      {\vartheta\atop1}
      \right)
  $$
  and the linear Weingarten coefficients of \reqn{lw}
  of the parallel nets $({\fk f,t})(\vartheta)$ become
  $$
    \alpha(\vartheta)
    = \alpha+2\beta\vartheta+\gamma\vartheta^2,
      \enspace
    \beta(\vartheta)  = \beta+\gamma\vartheta
      \enspace{\rm and}\enspace
    \gamma(\vartheta) = \gamma.
  $$
\item{} Thus, in the generic case,
  $\vartheta=-{\beta\over\gamma}$ yields
  a constant Gauss curvature net,
  which has two parallel constant mean curvature nets at
  $\vartheta=-{\beta\over\gamma}
    \pm{1\over\gamma}\sqrt{\beta^2-\alpha\gamma}$
  as soon as the root is real:
  we recover the classical Bonnet theorem in this case,
  cf \href[ref.imdg]{[14, \S2.7.4]}.
\item{} If however $\gamma=0$
  then the parallel family consists of the parallel
  constant harmonic mean curvature nets of a minimal net
  at $\vartheta=-{\alpha\over2\beta}$,
  as discussed in \href[expl.chmc]{Example 4.5}.

\item{(3)} In the indefinite case of hyperbolic and de Sitter
  spaces we base the analysis on an orthonormal basis again,
  so that
  $$
    B = \left(
      {\cosh\vartheta\atop\sinh\vartheta}\,
      {\sinh\vartheta\atop\cosh\vartheta}
      \right)
  $$
  now parametrizes the parallel family.
  To obtain a convenient representation of the linear
  Weingarten coefficients of \reqn{lw} for the parallel
  nets $({\fk f,t})(\vartheta)$ we distinguish three cases.
\item[4em]{(i)} If $|{\alpha+\gamma\over 2}|>|\beta|$
  we write
  ${\alpha-\gamma\over 2}=\mu$,
  ${\alpha+\gamma\over 2}=\varrho\cosh2\omega$ and
  $\beta=\varrho\sinh2\omega$
  to obtain
  $$
    \alpha(\vartheta) = \mu+\varrho\cosh2(\vartheta+\omega),
      \enspace
    \beta(\vartheta)  = \varrho\sinh2(\vartheta+\omega)
      \enspace{\rm and}\enspace
    \gamma(\vartheta) = -\mu+\varrho\cosh2(\vartheta+\omega).
  $$
  Thus $\vartheta=-\omega$ yields one constant Gauss curvature
  net, which has two parallel constant mean or constant
  harmonic mean curvature nets if
  $\varrho^2-\mu^2=\alpha\gamma-\beta^2<0$.
\item[4em]{(ii)} If ${\alpha+\gamma\over 2}=\pm\beta$
  we find, with ${\alpha-\gamma\over 2}=\mu\neq0$,
  $$
    \alpha(\vartheta) = \mu\pm\beta e^{\pm2\vartheta},
      \enspace
    \beta(\vartheta)  = \beta e^{\pm2\vartheta}
      \enspace{\rm and}\enspace
    \gamma(\vartheta) = -\mu\pm\beta e^{\pm2\vartheta}.
  $$
  Thus, as long as $\beta\neq0$,
  the parallel family contains
  no constant Gauss curvature net, but
  either one constant mean
  or one constant harmonic mean curvature net,
  depending on the sign of
  ${\alpha-\gamma\over\alpha+\gamma}$.
  On the other hand, $\beta=0$
  yields a parallel family of flat fronts,
  see \href[expl.ff]{Example 4.3}.
\item[4em]{} Note that we obtain discrete analogues
  of {\it linear Weingarten surfaces of Bryant type\/}
  in hyperbolic or de Sitter space here:
  with $\varepsilon=-({\fk pp})$
  and $\kappa=-({\fk qq})=-\varepsilon$
  and assuming, without loss of generality,%
  \footnote{A change ${\fk p}\to-{\fk p}$ of orientation
    reverses the sign of $H$.}
  that ${\alpha+\gamma\over2}=-\beta$,
  the linear Weingarten condition reads
  $$
    0
    = (\mu-\beta)(K-1) + 2\beta\,(H-1)
    = (\mu-\beta)\varepsilon\,K_{int} + 2\beta\,(H-1).
  $$
  As in the smooth case, see \href[ref.bjr11]{[8, Sect 4]},
  a characteristic feature of these nets is that
  one of the linear sphere complexes ${\fk k}^\pm$
  consists of spheres touching the infinity sphere
  ${\fk p-q}$:
  $$
    W
    = (\mu-\beta)\,{\fk q\odot q}
    + 2\beta\,{\fk q\odot p}
    - (\mu+\beta)\,{\fk p\odot p}
    = {\fk k^+\odot k^-}
  $$
  with ${\fk k}^+={\fk p-q}$
  and ${\fk k}^-=(\beta-\mu)\,{\fk q}-(\beta+\mu)\,{\fk p}$.
  This observation leads to a geometric Weierstrass type
  representation for these linear Weingarten nets,
  cf \href[ref.hrsy09]{[15]}.
\item[4em]{(iii)} If $|{\alpha+\gamma\over 2}|<|\beta|$
  we write
  ${\alpha-\gamma\over 2}=\mu$,
  ${\alpha+\gamma\over 2}=\varrho\sinh2\omega$ and
  $\beta=\varrho\cosh2\omega$
  to find
  $$
    \alpha(\vartheta) = \mu+\varrho\sinh2(\vartheta+\omega),
      \enspace
    \beta(\vartheta)  = \varrho\cosh2(\vartheta+\omega)
      \enspace{\rm and}\enspace
    \gamma(\vartheta) = -\mu+\varrho\sinh2(\vartheta+\omega).
  $$
  Thus the family does not contain any nets of constant Gauss
  curvature but a constant mean and a constant harmonic mean
  curvature net are obtained when
  $\sinh2(\vartheta+\omega)=\mp{\mu\over\varrho}$.
  When $\mu=0$ these coincide and the net is minimal.

Thus in any family of parallel (non-tubular) linear Weingarten
nets occurs at least one of the special nets that the Lawson
transformation was discussed for in the above examples.
In particular, our analysis shows that the genericity issue
of the Lawson transformation only occurs in the cases of
parallel families of intrinsically flat surfaces
in non-zero ambient curvature.%
\footnote{Note that ``intrinsically flat nets'' in a flat
  ambient geometry would be tubular.}

\goodbreak\vskip 3ex
{\usefont[cmbx10 scaled 1200]References}
\vglue 1ex
{\frenchspacing

\item{\htag[ref.bl29]{1.}}
 W Blaschke:
 {\it Vorlesungen \"uber Differentialgeometrie III\/};
 Springer Grundlehren XXIX, Berlin (1929)

\item{\htag[ref.bosu08]{2.}}
 A Bobenko, Y Suris:
 {\it Discrete Differential Geometry: Integrable Structure\/};
 Graduate Studies in Mathematics 98, Amer Math Soc,
 Providence (2008)

\item{\htag[ref.bpw10]{3.}}
 A Bobenko, H Pottmann, J Wallner:
 {\it A curvature theory for discrete surfaces
   based on mesh parallelity\/};
 Math Ann 348, 1--24 (2010)

\item{\htag[ref.bjl11]{4.}}
 A Bobenko, U Hertrich-Jeromin, I Lukyanenko:
 {\it Discrete constant mean curvature nets
   via K\"onigs duality\/};
 Discr Comput Geom 52, 612--629 (2014)

\item{\htag[ref.bdpp11]{5.}}
 F Burstall, N Donaldson, F Pedit, U Pinkall:
 {\it Isothermic submanifolds of symmetric $R$-spaces\/};
 J reine angew Math 660, 191-243 (2011) 

\item{\htag[ref.bjrs08]{6.}}
 F Burstall, U Hertrich-Jeromin, W Rossman, S Santos:
 {\it Discrete surfaces of constant mean curvature\/};
 RIMS Kyokuroku Bessatsu 1880, 133--179 (2014)

\item{\htag[ref.bjr10]{7.}}
 F Burstall, U Hertrich-Jeromin, W Rossman:
 {\it Lie geometry of flat fronts in hyperbolic space\/};
 C R 348, 661--663 (2010)

\item{\htag[ref.bjr11]{8.}}
 F Burstall, U Hertrich-Jeromin, W Rossman:
 {\it Lie geometry of linear Weingarten surfaces\/};
 C~R 350, 413--416 (2012)

\item{\htag[ref.cecil]{9.}}
 T Cecil:
 {\it Lie sphere geometry. With applications to submanifolds\/};
 Universitext, Springer,  New York (1992)

\item{\htag[ref.de11a]{10.}}
 A Demoulin:
 {\it Sur les surfaces $R$ et les surfaces $\Omega$\/};
 C R 153, 590--593 \& 705--707 (1911)

\item{\htag[ref.de11b]{11.}}
 A Demoulin:
 {\it Sur les surfaces $\Omega$\/};
 C R 153, 927--929 (1911)

\item{\htag[ref.do06]{12.}}
 A Doliwa:
 {\it Generalized isothermic lattices\/};
 J Phys A: Math Theor 40, 12539-12561 (2007)

\item{\htag[ref.jtz97]{13.}}
 U Hertrich-Jeromin, E Tjaden, M Z\"urcher:
 {\it On Guichard's nets and cyclic systems\/};
 EPrint arXiv:dg-ga/9704003 (1997)

\item{\htag[ref.imdg]{14.}}
 U Hertrich-Jeromin:
 {\it Introduction to M\"obius differential geometry\/};
 London Math Soc Lect Notes Series 300,
  Cambridge Univ Press, Cambridge (2003)

\item{\htag[ref.hrsy09]{15.}}
 T Hoffmann, W Rossman, T Sasaki, M Yoshida:
 {\it Discrete flat surfaces and linear Weingarten surfaces
   in hyperbolic $3$-space\/};
 Trans Amer Math Soc 364, 5605--5644 (2012)

\item{\htag[ref.muni06]{16.}}
 E Musso, L Nicolodi:
 {\it Deformation and applicability of surfaces in
   Lie sphere geometry\/};
 T\^ohoku Math J 58, 161--187 (2006)

\item{\htag[ref.pate88]{17.}}
 R Palais, C-L Terng:
 {\it Critical point theory and submanifold geometry\/};
 LNM 1353, Springer, Berlin (1988)

}
\vfill
\bgroup
\small
\def\addwd{\hsize=.36\hsize}
\def\spacer{\vtop{\addwd~\vfill}}
\def\udo{\vtop{\addwd
  U Hertrich-Jeromin\\
  Technische Universit\"at Wien\\
  Wiedner Hauptstra\ss{}e 8-10/104\\
  1040 Wien (Austria)\\
  Email: udo.hertrich-jeromin@tuwien.ac.at
  }}
\def\fran{\vtop{\addwd
  F Burstall\\
  Department of Mathematical Sciences\\
  University of Bath\\
  Bath, BA2 7AY (United Kingdom)\\
  Email: f.e.burstall@bath.ac.uk
  }}
\def\wayne{\vtop{\addwd
  W Rossman\\
  Department of Mathematics\\
  Kobe University\\
  Rokko, Kobe 657-8501 (Japan)\\
  Email: wayne@math.kobe-u.ac.jp
  }}
\hbox to \hsize{\hfil \fran \hfil \udo \hfil}\vskip 4ex
\hbox to \hsize{\hfil \wayne \hfil \spacer \hfil}
\egroup

\bye